\documentclass[reqno,10pt]{amsart} 
\usepackage{amsmath, amsfonts, amsthm, amssymb, verbatim,dsfont,stmaryrd}
\usepackage[mathcal]{euscript}
\usepackage{epsfig}

\newtheorem{definition}{Definition}[section]
\newtheorem{lemma}[definition]{Lemma}
\newtheorem{theorem}[definition]{Theorem}
\newtheorem{proposition}[definition]{Proposition}

\numberwithin{equation}{section}

\newcommand \Curl {\text{Curl}}
\newcommand \barC {\underline C}

\newcommand \Cbar {\overline C}

\newcommand \rhot {\widetilde \rho}

\newcommand \lbrac \llbracket 
\newcommand \rbrac \rrbracket 
\newcommand \R    {\mathbb{R}}

\newcommand \RR {\mathbb{R}}
 
\newcommand \del  {{\partial}}  

\newcommand \be {\begin{equation}}
\newcommand \ee {\end{equation}}

\newcommand \Fcal {\mathcal{F}}

\newcommand \eps {\epsilon}

\newcommand \la \langle
\newcommand \ra \rangle 
  
\let\oldmarginpar\marginpar
\renewcommand\marginpar[1]{\-\oldmarginpar[\raggedleft\footnotesize #1]%
{\raggedright\footnotesize #1}}
\newcommand \auth {\textsc}

\begin{document} 

\bibliographystyle{plain}
\title[A priori estimates for relativistic fluid flows with vacuum]
{Lagrangian formulation and a priori estimates for relativistic fluid flows with vacuum} 
\author[J. Jang, P.G. L{\tiny e}Floch, and N. Masmoudi]{
Juhi Jang$^1$, Philippe G. L{\smaller e}Floch$^2$, \and 
Nader Masmoudi$^3$ 
}
\thanks{
$^1$~Department of Mathematics, University of  Southern California, Los Angeles CA 90058, USA. Email: juhijang@usc.edu. 
\\
$^2$~Laboratoire Jacques-Louis Lions \& Centre National de la Recherche Scientifique, 
Universit\'e Pierre et Marie Curie, 4 Place Jussieu, 75252 Paris, France.
Email: contact@philippelefloch.org.
\\
$^3$~Courant Institute of Mathematical Sciences, New York University, 251 Mercer Street, NY 10011 New York, USA. Email: masmoudi@cims.nyu.edu. 
\\ 
2000\textit{\ AMS Subject Classification.} Primary: 35L65. Secondary: 76L05, 76N. \textit{Keywords. Relativistic fluid, vacuum state, free boundary, Lagrangian formulation, weighted energy}}
\begin{abstract}  We study the evolution of a compressible fluid surrounded by vacuum and introduce a new symmetrization in Lagrangian coordinates that allows us to encompass both relativistic and non-relativistic fluid flows. The problem under consideration is a free boundary problem of central interest in compressible fluid dynamics and, from the mathematical standpoint, the main challenge to be overcome lies in the loss of regularity in the fluid variables near the free boundary.  Based on our Lagrangian formulation, we establish the necessary a priori estimates in weighted Sobolev spaces which are adapted to this loss of regularity.
\end{abstract}
\maketitle  



\section{Introduction} 

We study here the Euler equations describing the evolution of a relativistic compressible fluid, that is, the system 
\be
\label{10} 
\aligned
& \del_t \big( \rhot - \eps^2 p(\rho) \big) + \text{div} (\rhot \, u) = 0, 
\\
& \del_t \big( \rhot \, u \big) + \text{div} \big( \rhot \, u \otimes u) + \text{grad} \big( p(\rho) \big) = 0, 
\endaligned
\ee 
in which the mass density $\rho=\rho(t,x)$ and
the velocity vector of the fluid $u=u(t,x)$ (with $t \geq 0$ and $x \in \RR^3$)  are the main unknowns and 
satisfy the physical constraints 
\be
\label{11} 
\rho \geq 0, \qquad |u| < 1/\eps. 
\ee
The parameter $1/\eps$ represents the light speed and, in \eqref{10}, the pressure $p=p(\rho)$ is a given function of the density, while the ``modified density'' $\rhot$ is defined by  
\be
\label{12} 
\rhot=\rhot(\rho, u) := {\rho + \eps^2 p \over 1 - \eps^2 \, |u|^2}.  
\ee 
We also set $x= (x_ i)_{1 \leq i \leq 3}$ and use the standard notation for the divergence 
$\text{div} u := \sum_ i \del_{x_ i} u_ i = \del_i u_i$
(with implicit summation on $i$) and for the gradient $\text{grad}( p ) = \big( \del_i p \big)_{1 \leq i \leq 3}$.

Under the standard physical assumption that $p'(\rho) \geq 0$ (and vanishes if and only if $\rho = 0$), the Euler equations \eqref{10} away from the vacuum state form a strictly hyperbolic system of four conservation laws, which, however, is {\sl non-strictly hyperbolic} at the vacuum $\rho = 0$.  
 We are interested in the evolution of  a 
compressible fluid region surrounded by vacuum, in particular when a fluid is continuously in contact with vacuum. This is a classical problem in fluid dynamics
and, from the mathematical standpoint, the main technical challenge to be overcome lies in the loss of regularity in the fluid variables near the  free boundary between the fluid and the vacuum region. Specifically, we require that the normal acceleration of the fluid near the boundary is non-vanishing and bounded: 
\be\label{physical_vacuum}
\barC \leq | \del_\nu p'(\rho) | \leq \Cbar 
\ee
for some constants $0 < \barC \leq \Cbar < +\infty$, where $\nu \in \RR^3$ denotes the normal unit vector to the 
free fluid-vacuum boundary. This vacuum boundary condition, the so-called ``physical vacuum'' boundary, can be realized by some self-similar solutions and stationary solutions for different physical systems such as Euler equations with damping and Euler-Poisson systems for gaseous stars \cite{J12,JM0,LY2,Makino2012,Y}. 

Let us mention several earlier works on the above problem which attracted a lot of attention in recent years. 
Coutand and Shkoller \cite{CS1,CS2} successfully established an existence result for non-relativistic compressible fluids  (that is, the system \eqref{10} with $\eps =0$) by degenerate parabolic regularizations, while, independently, Jang and Masmoudi developed a hyperbolic-type weighted energy estimates for all spatial derivatives including normal derivatives in order to prove the existence of solutions in one space dimension \cite{JM1} as well as in several space dimensions \cite{JM2}. We also mention that in a recent work \cite{Makino2012}, Makino addressed some existence result for the Euler-Poisson system based on the Nash-Moser-Hamilton theory. We refer to \cite{JM0,JM2} for a historical background and a bibliography on the subject.  

As far as relativistic fluids are concerned, earlier investigations on compactly supported solutions to the relativistic Euler equations include Makino and Ukai \cite{MU} and, for the equations in several space variables,  LeFloch and Ukai \cite{LU}. In these works, a stronger regularity property is implied on the fluid variables near the free boundary between the fluid and the vacuum. A general existence theory for relativistic compressible fluids encompassing the above vacuum condition \eqref{physical_vacuum} is therefore still lacking. 

The goal of this article is to present a new Lagrangian formulation of the relativistic Euler equations and is to derive the necessary a priori bounds satisfied by solutions subject to \eqref{physical_vacuum} based on such our formulation. 

An outline of the paper is as follows. In Section \ref{sec:2}, we present the compressible fluid flow equations, discuss its reduction to a second--order hyperbolic system in Lagrangian variables, and also derive the relativistic vorticity equation. In Section~\ref{sec:3}, we introduce 
the free boundary problem of interest and present the a priori estimates. Section ~\ref{sec:4} includes discussion on non-relativistic flows as well as the existence theory for special cases such as radially symmetric flows.  


\section{Lagrangian formulation for relativistic fluid flows}
\label{sec:2}

\subsection{Equations of state}
 
The Euler equations in Minkowski spacetime read 
\be
\label{Euler}
\aligned 
&\del_t\Bigl( {\rho + \eps^2 p \over 1 - \eps^2 \, |u|^2} - \eps^2 \, p \Bigr)
  + \del_k \Bigl( {\rho + \eps^2 p \over 1 - \eps^2 \, |u|^2} \, u_k \Bigr) = 0,
\\
& \del_t \Bigl( {\rho + \eps^2 p \over 1 - \eps^2 \, |u|^2} \, u_j \Bigr) 
  + \del_k \Bigl( {\rho + \eps^2 p \over 1 - \eps^2 \, |u|^2} \, u_j u_k + p \, \delta_{jk} \Bigr) = 0,  
\endaligned  
\ee
where  $(\rho,u): [0,T] \times \RR^3 \to \RR_+ \times \RR^3$ is the main unknown, defined on some time interval $[0,T)$. 
As pointed out in the introduction, it is convenient to introduce the variable $\rhot$,  so that the Euler equations read 
\be
\label{Euler28}
\aligned
&\del_t \bigl(\rhot  - \eps^2 \, p \bigr)
  + \del_k \bigl( \rhot  \, u_k \bigr) = 0,  
\\
& \del_t \bigl( \rhot  \, u_j \bigr) 
  + \del_k \bigl( \rhot \, u_j u_k + p \, \delta_{jk} \big) = 0.  
\endaligned
\ee
Observe that, by letting formally $\eps \to 0$ in \eqref{Euler}
we find $\rhot \to \rho$ and we recover the non-relativistic Euler equations.  

As it is required by the physics of the problem, the {\sl sound speed} $c(\rho):= p'(\rho)^{1/2}$ is assumed to be real and smaller than the light speed, that is, 
\be
\label{hyperbolic} 
0 < c(\rho) < \eps^{-1}  \quad \text{ provided } \rho >0. 
\ee 
For concreteness,  the pressure $p$ is assumed to be a power-law of the particle number $N$, that is, 
$$
p = {a \over \gamma - 1} N^\gamma, 
$$
where $N$ is related to the energy density $$\rho = N +  \eps^2 {a \over \gamma - 1} N^\gamma$$ and 
$\gamma \in (1,2)$ is a constant referred to as the adiabatic exponent of gases and $a>0$ is a normalization constant. 
Hence, the pressure is determined implicitly by 
$$
\rho = \eps^2 p + \kappa \, p^{1 / \gamma}, \qquad \kappa^\gamma := {\gamma - 1 \over a}.  
$$
It is easily checked that the sound speed does not exceed the light speed,  
since
$$
c^2 = p'(\rho) = {1 \over { \eps^2 } + { \gamma-1 \over a \gamma}N^{1-\gamma } } \leq {1 \over \eps^2}.   
$$  

We also define the function $$h=h(\rho) \;\text{ by }\; dh := {1 \over N} dp$$ and, from now on, adopt the normalization $a:= \gamma - 1$, so that the 
equation of state of the fluid finally reads 
\be
\label{eq:pressure}
p(\rho) = N^\gamma, \qquad h(\rho) = {\gamma \over \gamma - 1} N^{\gamma-1}, 
\qquad 
\text{ with } 
\rho = N + \eps^2N^\gamma, 
\ee
where $\gamma \in (1,2)$. 


\subsection{The energy equation and the number density equation}

We recall the energy pair $(V,H)$ associated with the relativistic Euler equations
\be
\label{energy}
\del_t V + \del_j H_j = 0, 
\ee
with 
$$
\aligned
& V := \eps^{-2} \left( (1 + \kappa\eps^2) (\rhot - \eps^2 p) 
               - {N(\rho) \over (1 - \eps^2 |u|^2)^{1/2}} \right),
\\
& H_j :=  \eps^{-2} \left(  (1 + \kappa\eps^2) \rhot  \, u_j
             - {N(\rho) u_j  \over (1 - \eps^2 |u|^2)^{1/2}}\right),
\endaligned 
$$
where 
$$
N(\rho) := \exp\left(\int_1^\rho {ds \over s +\eps^2 p(s)}\right), 
\qquad 
\kappa := \int_0^1 {p(s) \over s^2} \, ds. 
$$ 
Here the function $N$ is determined so that, as $\eps \to 0$, the pair $(V,H)$
tends to the standard 
energy pair of the non-relativistic Euler equations. Indeed, as $\eps \to 0$ we have 
$$
\aligned 
& V \sim  {1 \over 2} \rho |u|^2 + \rho \int_0^\rho {p(s) \over s^2} \, ds, 
\\
& H_j \sim u_j \, \big( \rho \, V + p(\rho) \big). 
\endaligned 
$$
It can be checked that 
the above function is strictly convex \cite{MU} in the conservative variable $\omega=(\rho,\rho u)$, and 
$$
\nabla_\omega^2 V  \geq C_1 \qquad    \text{ away from the vacuum,} 
$$
where the constant $C_1$ is uniform on every compact subset of $\big\{ \rho > 0, \, |u| < \eps \big\}$, 
excluding therefore the vacuum.

From the mass density equation in \eqref{Euler} and the energy equation \eqref{energy},  we deduce  that any solution $(\rho,u): \RR_+ \times \RR^3 \to \RR_+ \times \RR^3$ to the Euler equations
\eqref{Euler} also satisfies the following number density equation
\be
\label{energy_0}
\aligned
& \del_t g + \del_j \big( g \, u_j \big) = 0, 
\\
& g := N \Gamma, \qquad \Gamma =\Gamma(u) := (1 - \eps^2 |u|^2)^{-1/2}. 
\endaligned
\ee  
In the following we will work with \eqref{energy_0} together the second equation in \eqref{Euler28}. Note that the Cauchy problem 
is posed by prescribing, at the initial time $t=0$, the initial density $\rho_0$ and the initial velocity $u_0$ of the fluid
\be
\label{initial0}
\rho(0,x) = \rho_0(x), \qquad u(0,x) = u_0(x), \qquad x \in \RR^3
\ee
with, of course, 
\be
\label{initial1} 
\rho_0 \geq 0, \qquad |u_0| < 1/\eps. 
\ee
We are interested in the situation where the density is positive in some smooth open set $\Omega \subset \RR^3$ and vanishes identically outside this set, i.e. 
\be
\label{initial3} 
\rho_0 \, \begin{cases}
>0,   & x \in \Omega, 
\\
= 0, & x \in \RR^3 \setminus \Omega. 
\end{cases}
\ee


\subsection{Lagrangian coordinates and notation}

We are going to now reformulate the fluid equations above in terms of the {\sl Lagrangian coordinates} $\eta_j= \eta_j(t,x)$
 defined by the following ordinary differential equation with prescribed initial data: 
\be
\label{vvv}
\aligned
& \del_t \eta_j := u_j (t, \eta), 
\\
& \eta_j(0,x) := x_j. 
\endaligned
\ee
We introduce the Jacobian matrix of this transformation, that is, $A^{-1} := \big( \nabla_x \eta \big) = \big( \del_i  \eta_j \big)$ and $J := \det \big( \nabla_x \eta \big)$. We use Einstein's summation convention and the notation $F,_k$ to denote the $k$-th partial derivative of $F$: $\partial_kF$. Both expressions will be used throughout the paper.  Differentiating the inverse of deformation tensor,  since $A\cdot
[D\eta] =I$, one obtains  
\begin{equation}\label{DA}
\partial_t A^k_i=- A^k_r\partial_t \eta^r,_sA^s_i\;;\quad \partial_lA^k_i=-
A^k_r\partial_l\eta^r,_sA^s_i\,. 
\end{equation} 
Differentiating the Jacobian determinant, one obtains 
\begin{equation}\label{DJ}
\partial_tJ = JA^s_r \partial_t \eta^r,_s\;; \quad \partial_l J= JA^s_r \partial_l\eta^r,_{s}\,.
\end{equation}
For the cofactor matrix $JA$, from \eqref{DA} and \eqref{DJ},  one obtains the following Piola identity: 
\begin{equation}\label{Piola}
(JA^k_i),_k=0\,.
\end{equation}

For a given vector field $F$, we use $DF$, $\text{div}F$, $\text{curl}F$ to denote its full gradient, its divergence, and its curl: 
\[
\begin{split}
[DF]^i_j & \equiv F^i,_j \\ 
\text{div}F \;&\equiv  F^r,_r\\
[\text{curl}F]^i & \equiv \epsilon_{ijk} F^k,_j
\end{split}
\] 
where $\epsilon_{ijk}$ is the Levi-Civita symbol: it is 1 if $(i,j,k)$ is an even permutation of $(1,2,3)$, -1 if $(i,j,k)$ is an odd permutation of $(1,2,3)$, and 0 if any index is repeated. 

We introduce the following Lie derivatives along the flow map $\eta$:
\[
\begin{split}
 [D_\eta F]^i_r&\equiv A^s_rF^i,_s\\ \text{div}_\eta F\;&\equiv A^s_rF^r,_s \\
[\text{curl}_\eta F]^i &\equiv \epsilon_{ijk}A^s_jF^k,_s
\end{split}
\]
which indeed correspond to Eulerian full gradient, Eulerian divergence, and Eulerian curl written in Lagrangian coordinates. 
In addition, it is convenient to introduce the anti-symmetric  curl matrix $\text{Curl}_\eta F$: 
\[
[\text{Curl}_\eta F]^i_j \equiv  A^s_jF^i,_s - A^s_iF^j,_s \,.
\]
Note that $\text{Curl}_\eta F$ is a matrix version of a vector $\text{curl}_\eta F$ and that  $|\text{Curl}_\eta F|^2 
=2 |\text{curl}_\eta F|^2$ holds. We will use both $\text{curl}_\eta $ and $\text{Curl}_\eta $.  We end this section by recalling the following property of the Lagrangian curl: 
$$\textit{if} \;\; \omega^k= A^r_k f,_r,  \;\;  \text{curl}_\eta \omega =0.$$

\subsection{Lagrangian formulation}

By introducing the {\sl modified velocity} 
$$
\chi^j := (1 + \eps^2 h) \, \Gamma \del_t \eta^j,
$$
we arrive at the following equivalent formulation of the Euler equations 
\be
\label{290}
\aligned
\del_t g + g A^k_j \del_k u_j = 0, 
\qquad 
g \, \del_t \chi^j + A_j^k \del_k p = 0.
\endaligned
\ee

Furthermore, from  $\del_t J - J A_j^k \del_k u_j = 0$,
we deduce that $g_0 := g J = \Gamma N J$. Therefore, $N$ can be expressed in terms of the main unknown $\eta$ by 
$$
N = {g_0 \over \Gamma J} =  {g_0 \over J \, (1 - \eps^2 |\del_t \eta|^2)^{-1/2}}. 
$$
In turn, the second equation in \eqref{290} reads 
\be
\label{300} 
g_0 \del_t \Big(  (1 + \eps^2 h) \, \Gamma \del_t \eta^j \Big) + J A_j^k \del_k N^\gamma = 0,
\ee
in which the spatial derivative terms take also the form
$$
 J A_j^k \del_k N^\gamma 
= \del_k \Big( g_0^\gamma A_j^k J^{1 - \gamma} \Gamma^{-\gamma} \Big) 
$$
and, using the expression of $\Gamma$, 
$$
 J A_j^k \del_k N^\gamma 
= \del_k \Big( g_0^\gamma A_j^k J^{1 - \gamma} \Big)  \Gamma^{-\gamma}
 - g_0^\gamma A_j^k J^{1 - \gamma} \, \gamma \Gamma^{-\gamma-1} \eps^2 \Gamma^3 \del_t \eta^i \del_t \del_k \eta^i. 
$$
On the other hand, the time derivative in \eqref{300} takes the form
\be
\label{302} 
\aligned
& \del_t \Big(  (1 + \eps^2 h) \, \Gamma \del_t \eta^j \Big) 
\\
& = (1 + \eps^2 h) \, \Gamma \del_t^2 \eta^j 
 + 
\del_t \eta^j \Big( (1+ (2-\gamma) \eps^2 h) \, \eps^2 \Gamma^3 \del_t \eta^i \del_t^2 \eta^i - (\gamma-1) \eps^2 h \Gamma \del_t \log J \Big)  
\\
& = \Big( (1 + \eps^2 h) \,  \delta_i^j 
           +  (1+ (2-\gamma) \eps^2 h) \, \eps^2 \Gamma^2 \del_t \eta^i \del_t \eta^j \Big) \, \Gamma \del_t^2 \eta^i
- (\gamma-1) \eps^2 h \Gamma  \del_t \eta^j  \del_t \log J. 
\endaligned
\ee
The latter term can be expressed as a second-order term in $\eta$, by writing 
$
\del_t J = J A_i^k \del_k \del_t \eta^i 
$ so that 
\be
\label{303}
g_0 (\gamma-1) \eps^2 h \Gamma  \del_t \eta^j  \del_t \log J
= 
\gamma g_0^\gamma \Gamma^{- \gamma} J^{1-\gamma} A_i^k \del_k \del_t \eta^i \del_t \eta^j. 
\ee

Plugging \eqref{302}--\eqref{303} in the equation \eqref{300} and collecting the terms, 
we thus find a second-order equation in $\eta$ 
$$
g_0 B_i^j \del_t^2 \eta^i + g_0^\gamma C_{ij}^k \del_k \del_t \eta^i + \del_k \big( g_0^\gamma A_j^k J^{1-\gamma} \big) =0, 
$$
in which the coefficients are given by 
\[
\aligned
B_i^j := & \Big( (1 + \eps^2 h) \,  \delta_i^j 
           +  (1+ (2-\gamma) \eps^2 h) \, \eps^2 \Gamma^2 \del_t \eta^i \del_t \eta^j \Big) \, \Gamma^{\gamma+1}, 
\\
C_{ij}^k := & - \gamma \eps^2 \Gamma^2 J^{1-\gamma} \Big( A_i^k \del_t \eta^j + A_j^k \del_t \eta^i \big). 
\endaligned
\]
Finally, by letting 
$$
g_0 = w^\alpha, \qquad g_0^\gamma = w^{1+\alpha} \quad \text{where}\quad\alpha = {(\gamma-1)}^{-1}
$$
we arrive at the following {\sl second--order formulation in Lagrangian coordinates}
\be
\label{3010}
w^\alpha \, B_i^j \del_t^2 \eta^i + w^{1+\alpha} C_{ij}^k \del_k \del_t \eta^i + \del_k \big( w^{1+\alpha} A_j^k J^{-1/\alpha} \big) =0.  
\ee
Importantly, we have the symmetry property $B_i^j = B_j^i$ and $C_{ij}^k = C_{ji}^k$, while $B_i^j$ is positive definite. 

In Lagrangian coordinates, we prescribe the reference density function $\rho_0 \geq 0$
(which determines the initial data $g_0= w^\alpha$) so that it is positive in some smooth open set $\Omega \subset \RR^3$ and vanishes identically outside this set (cf.~\eqref{initial3} above)  
and we can then pose the Cauchy problem of interest by requiring that 
\be
\label{initial000}
\eta(0,x) = x, \qquad \eta_t(0,x) = \eta_1(x), \qquad x \in \Omega
\ee
for some data $\eta_1$ (which is precisely the velocity data $u_0$ in \eqref{initial000} expressed in Lagrangian coordinates).


\subsection{Relativistic vorticity} 

One additional set of equations will be required in our analysis. Observe that the second equation in \eqref{290} can be rewritten as 
\be
\label{chi-h}
\Gamma \del_t\chi^j +A^k_j\del_k h=0. 
\ee
Note that this equation allows us to control the spatial derivatives $\del_k h$ by the time derivative $\del_t \chi$. 
By taking the curl of that equation in Lagrangian coordinates, we obtain $\text{Curl}_\eta (\Gamma\del_t\chi)=0$, implying 
$$
\Gamma\text{Curl}_\eta\del_t\chi+ [\text{Curl}_\eta,\Gamma]\del_t\chi=0, 
$$
with 
$$ 
[\text{Curl}_\eta,\Gamma]\del_t\chi =A^l_i\Gamma,_l\del_t\chi^j - A^l_j\Gamma,_l \del_t\chi^i.
$$ 
Since $[\del_t,\text{Curl}_\eta]\chi=\del_t A^l_i\chi^j,_l -\del_t A^l_j\chi^i,_l$, this equation can be written as 
\be
\label{LV}
\del_t\text{Curl}_\eta\chi=[\del_t,\text{Curl}_\eta]\chi-\Gamma^{-1} [\text{Curl}_\eta,\Gamma]\del_t\chi, 
\ee
 which we refer to as the {\sl Lagrangian relativistic vorticity equation.}

By integrating \eqref{LV} in time, we deduce that 
\be
\label{curl-chi}
\text{Curl}_\eta\chi = \text{Curl}_\eta\chi \big|_{t=0}+\int_0^t  [\del_t,\text{Curl}_\eta]\chi\, ds-  \int_0^t  \Gamma^{-1}  [\text{Curl}_\eta,\Gamma]\del_t\chi \,ds.
\ee
For the purpose of the energy estimates, we will need to derive the equation for the curl of $\eta$ and estimate them. From the definition of $\chi$, we see that 
\[
\text{Curl}_\eta\chi= \Gamma(1+\eps^2 h)\text{Curl}_\eta\del_t\eta +{[\text{Curl}_\eta,\Gamma(1+\eps^2h)]\del_t\eta}. 
\]
Hence, \eqref{curl-chi}  reads as 
\begin{equation}
\begin{split}\label{curl-eta}
&\text{Curl}_\eta\del_t\eta+\boxed{[ \Gamma(1+\eps^2 h)]^{-1}[\text{Curl}_\eta,\Gamma(1+\eps^2h)]\del_t\eta}
\\
&= [ \Gamma(1+\eps^2 h)]^{-1} \Big(\text{Curl}_\eta\chi \big|_{t=0} 
+\int_0^t  [\del_t,\text{Curl}_\eta]\chi\, ds- \int_0^t  \Gamma^{-1}  [\text{Curl}_\eta,\Gamma]\del_t\chi \,ds\Big). 
\end{split}
\end{equation}
We observe that the boxed term in \eqref{curl-eta}  is {\sl not} of lower order. By using \eqref{chi-h}, we rewrite it so that it does not contain two spatial derivatives of $\eta$: 
\[
\begin{split}
\boxed{\text{ \ }^{\text{ \ } }} &=[ \Gamma(1+\eps^2 h)]^{-1}(\Gamma,_l (1+\eps^2 h)+\Gamma \eps^2 h,_l )(A^l_i\del_t\eta^j - A^l_j \del_t \eta^i)
\\
&=\eps^2\Gamma^2\del_t\eta^m\del_t\eta^m,_l \big( A^l_i\del_t\eta^j - A^l_j \del_t \eta^i \big) 
-\eps^2\Gamma(1+\eps^2 h)^{-1}(\del_t\chi^i\del_t\eta^j -\del_t\chi^j\del_t\eta^i)\\
&=\eps^2\Gamma^2( {\del_t\eta^m\del_t\eta^m,_l \big( A^l_i\del_t\eta^j} - A^l_j \del_t \eta^i \big) -\del_t^2\eta^i\del_t\eta^j+\del_t^2\eta^j\del_t\eta^i ).
\end{split}
\]
By rearranging terms, we write the curl equation  \eqref{curl-eta} as 

\begin{equation}
\begin{split}\label{curl-eta-new0}
&  [\text{D}_\eta  \partial_t \eta]_i^m  (\delta^j_m+   \eps^2\Gamma^2\del_t\eta^j\del_t\eta^m  )
-  (\delta^m_i+   \eps^2\Gamma^2\del_t\eta^i\del_t\eta^m  )  [\text{D}_\eta  \partial_t  \eta]_j^m  
\\
& + \eps^2  \Gamma^2 
  (\del_t^2\eta^j\del_t\eta^i  -\del_t^2\eta^i\del_t\eta^j )  \\  
&= [ \Gamma(1+\eps^2 h)]^{-1} \Big[\text{Curl}_\eta\chi \big|_{t=0} 
+\int_0^t  [\del_t,\text{Curl}_\eta]\chi\, ds- \int_0^t  \Gamma^{-1}  [\text{Curl}_\eta,\Gamma]\del_t\chi \,ds\Big]^j_i. 
\end{split}
\end{equation}

We next define  the symmetric matrix 
\[
S_{m}^j :=  (\delta^j_m+   \eps^2\Gamma^2\del_t\eta^j\del_t\eta^m  ) 
\]
and the anti-symmetric matrices  

\[
\begin{split}
R_{i}^j &:=\eps^2  \Gamma^2 
  (\del_t^2\eta^j\del_t\eta^i  -\del_t^2\eta^i\del_t\eta^j ) \\
X_{i}^j &:=  [ \Gamma(1+\eps^2 h)]^{-1} \Big[\text{Curl}_\eta\chi \big|_{t=0} 
+\int_0^t  [\del_t,\text{Curl}_\eta]\chi\, ds- \int_0^t  \Gamma^{-1}  [\text{Curl}_\eta,\Gamma]\del_t\chi \,ds\Big]_i^j 
\end{split}
\]
 Then, \eqref{curl-eta-new0} can be written as 
 \begin{equation}
\begin{split}\label{curl-eta-new}
&  [\text{D}_\eta  \partial_t   \eta]^m_i S_m^j - 
  S^m_i  [\text{D}_\eta  \partial_t   \eta]^m_j 
+ R_i^j = X_i^j 
\end{split}
\end{equation}

Notice that $S$ is symmetric and positive definite, hence, letting $U:= S^{-1}$, we get the following equivalent expression to the curl equation \eqref{curl-eta-new}: 
\[
 U^m_i [\text{D}_\eta  \partial_t   \eta]_m^j - 
   [\text{D}_\eta  \partial_t   \eta]^i_m U^j_m  
+ U_i^m R_m^l U_l^j =U_i^m X_m^l U_l^j. 
\]


\subsection{Relativistic Euler equations as a second-order hyperbolic system}  

So far, we have reformulated the relativistic Euler equations as a second-order quasi-linear hyperbolic system in Lagrangian coordinates, where $\eta= (\eta^j(t,x)) \in \R^3$ is the main unknown, and have identified the corresponding curl structure. We summarize such formulations in the following proposition. 

\begin{proposition} Suppose $(\rho, u)$ are smooth solutions to relativistic Euler equations \eqref{Euler} written in  Eulerian coordinates. Let $w^\alpha=N_0\Gamma_0$, where $N_0=N_0(\rho_0)$ is the initial particle number density determined by \eqref{eq:pressure} and $\Gamma_0=(1-\eps^2|u_0|^2)^{-1/2}$. Then the solution $\eta$ to the ODE \eqref{vvv} satisfies the following second-order quasi-linear hyperbolic system 
\be
\label{syste-0}
\aligned
& w^\alpha \, B_i^j \del_t^2 \eta^i + w^{\alpha+1} C_{ij}^k \del_k \del_t \eta^i + \del_k \big( w^{\alpha+1} A_j^k J^{-1/\alpha} \big) =0, 
\endaligned
\ee
where 
  \be\label{BC}
\aligned
B_i^j = & \Big( (1 + \eps^2 h) \,  \delta_i^j 
           +  (1+ (1-\tfrac1\alpha) \eps^2 h) \, \eps^2 \Gamma^2 \del_t \eta^i \del_t \eta^j \Big) \, \Gamma^{2+1/\alpha}, 
\\
C_{ij}^k = & - (1+\tfrac1\alpha) \eps^2 \Gamma^2 J^{-1/\alpha} \Big( A_i^k \del_t \eta^j + A_j^k \del_t \eta^i \big), 
\endaligned
\ee
and furthermore, admits the following structure 
\begin{equation}\label{curl-eta-new1}
 U^m_i [\text{D}_\eta  \partial_t   \eta]_m^j - 
   [\text{D}_\eta  \partial_t   \eta]^i_m U^j_m  
+ U_i^m R_m^l U_l^j =U_i^m X_m^l U_l^j,  
\end{equation} 
where 
  \be\label{URX}
\aligned
U^j_i &= (S^{-1})^j_i\quad \text{where }\;\; S_{m}^j =  (\delta^j_m+   \eps^2\Gamma^2\del_t\eta^j\del_t\eta^m  ), 
\\
R_{i}^j &=\eps^2  \Gamma^2 
  (\del_t^2\eta^j\del_t\eta^i  -\del_t^2\eta^i\del_t\eta^j ), \\
X_{i}^j &=  [ \Gamma(1+\eps^2 h)]^{-1} \Big[\text{Curl}_\eta\chi \big|_{t=0} 
+\int_0^t  [\del_t,\text{Curl}_\eta]\chi\, ds- \int_0^t  \Gamma^{-1}  [\text{Curl}_\eta,\Gamma]\del_t\chi \,ds\Big]_i^j. 
\endaligned
\ee
Here we recall that 
\be\label{recall}
 \chi^j = (1 + \eps^2 h) \, \Gamma \del_t \eta^j, \;\; h=(1+\alpha)w (\Gamma J)^{-1/\alpha},\;\; \Gamma= (1 - \eps^2 |\partial_t\eta|^2)^{-1/2}. 
\ee
Conversely, if $(\eta, \eta_t)$ (with $J$ being bounded away from zero and above) are smooth solutions to the above system, $(\rho, u)$ is a solution to the Eulerian system. 
\end{proposition}

We observe that this proposition can be justified at least away from vacuum, where smooth solutions are available in Eulerian coordinates; for instance, see \cite{LU, MU}. 

In the next section, based on the above reformulation \eqref{syste-0}--\eqref{recall} of the relativistic Euler equations in Lagrangian coordinates, we will establish the a priori estimates for smooth solutions in the presence of a physical vacuum. 

\section{The free boundary problem for the relativistic Euler system}\label{sec:3}
 
\subsection{Main result}

In this section, we consider a vacuum free boundary problem for relativistic Euler equations in Lagrangian coordinates. We first prescribe a class of $w$: $w$ is  the prescribed function in $\Omega$ with smooth boundary $\del\Omega$ and it vanishes at the boundary like a distance function:  
\be
\label{syste-4}
\aligned
&w= 0 \quad \text{ on } \del \Omega, 
\\
& \barC \, d(x, \del \Omega) \leq w \leq \Cbar \, d(x, \del\Omega). 
\endaligned
\ee
The regularity for $w$ will be specified in the next subsection. (See \eqref{F-w}.)

We can pose the Cauchy problem of interest by requiring that 
\be
\label{initial0-0}
\eta(0,x) = \eta_0(x), \qquad \eta_t(0,x) = \eta_1(x), \qquad x \in \Omega
\ee
for some given data $\eta_0, \eta_1$. Note that due to degeneracy of $w$ we indeed do not need to impose the boundary condition on $\del\Omega$. We are interested in the {\sl free boundary value problem} associated with the noninear hyperbolic systems \eqref{syste-0}, that is, we search for solutions that are supported in a domain $\Omega$ with smooth boundary $\del \Omega$. It is a moving boundary value problem because $\eta_t$, the velocity of fluids, is not necessarily zero along the boundary, and the moving domain in Eulerian coordinates is given by $\Omega(t)=\eta(t)(\Omega)$. 

Observe that the condition imposed near the boundary is singular in nature and 
special care will be required to handle derivatives of $\eta$, especially in the direction normal to the boundary.

For simplicity of the presentation, we consider the case when 
  the initial domain is taken as $$\Omega=\mathbb{T}^2\times
(0,1),
$$ 
where $\mathbb{T}^2$ is a two-dimensional period box in $x_1,x_2$. The result can be extended to the general case   in the same way as done in \cite{JM2}. 
  The initial boundary is given as 
$$
\partial\Omega=\{x_3=0\}\cup\{x_3=1\} \text{ as the reference vacuum boundary.}
$$
We use Latin letters  $i,j,k, \ldots$ to denote $1, 2, 3$
     and that we 
  use Greek letters $\beta,\kappa,\sigma, \tau$ to denote $1, 2$, only.
 We use $\partial_\tau^m$ to denote $\partial_1^{m_1}\partial_2^{m_2}$ and $|m|$ to denote $|m|=m_1+m_2$.  To any sufficiently regular function $\eta$ defined on $[0,T] \times \Omega$, we associate the following energy functionals (defined for any integer $N \geq 0$):   
$$
\aligned
E_N^{(I)}:=& \sum_{|m| + n \leq N} \int_\Omega  
 w^{\alpha+n} \del_\tau^m \del_3^n \eta_t^j\, B^j_i \,\del_\tau^m \del_3^n \eta_t ^i
 \, dx=:  \sum_{|m| + n \leq N} \mathcal{E}_{m,n}^{(I)},  
\\
E_N^{(II)} :=& \sum_{|m| + n \leq N} \int_\Omega  w^{\alpha+n+1} J^{-1/\alpha} |\text{div}_\eta \del_\tau^m \del_3^n \eta |^2 \, dx=:  \sum_{|m| + n \leq N} \mathcal{E}_{m,n}^{(II)},
\\
E_N^{(III)} :=& \sum_{|m| + n \leq N} \int_\Omega w^{\alpha+n+1}  \,  [D_\eta \del_\tau^m \del_3^n\eta]^j_m  {U}^i_m [D_\eta \del_\tau^m \del_3^n\eta]^j_i \,  dx =:  \sum_{|m| + n \leq N} \mathcal{E}_{m,n}^{(III)},
\\
E_N^{(IV)} :=& \sum_{|m| + n \leq N} \int_\Omega w^{\alpha+n+1}  \,  | \del_\tau^m \del_3^n \Curl_\eta \chi |^2 \,  dx. 
\endaligned
$$ 
Note that $E_N^{(II)}$ is bounded by $E_N^{(III)}$. 
 The total energy of interest is the sum 
$$
E_N= E_N^{(I)}+ E_N^{(III)} + E_N^{(IV)}.  
$$ 

Furthermore, the regularity of the weight function $w$ is determined by introducing the norms: 
\be\label{F-w}
\aligned
F_M[w] :=& \sum_{|m| + n \leq M} \int_\Omega w^{\alpha+n+1}  \,  |\del_\tau^m \del_3^n w |^2 \, dx,
\\
F_M^{(I)}[w] :=& \sum_{|m| + n \leq M} \int_\Omega  w^{\alpha+n+1}  \,  |D \del_\tau^m \del_3^n w |^2 \, dx. 
\endaligned
\ee

We now state the result on the a priori estimates for solutions of \eqref{syste-0} in the above energy spaces. 

\begin{theorem}[A priori estimates]
\label{theo1} Let $N \geq 2 \alpha + 9$ be fixed for given exponent $\alpha >0$ and let $w$ be given satisfying \eqref{syste-4} and 
$
F_N[ Dw] < \infty. 
$
Suppose $\eta$ and $\eta_t$ solve \eqref{syste-0} for $t\in[0,T]$ 
 with $E_N=E_N(\eta,\eta_t)<\infty$ and $1/C_0\leq J\leq C_0$ for some $C_0\geq 1$ for the initial data $\eta_0, \eta_1 :\RR^3 \to \RR^3$  in \eqref{initial0-0}
satisfying 
$
E_N[\eta_0, \eta_1] < \infty. 
$
  We further assume that $\eta$
and $\eta_t$ enjoy the a priori bound: for any $s=1,2,\text{ and }3$, 
\begin{equation}\label{Assumption}
 \sum_{|p|+q=0}^{[N/2]} |w^{q/2}\partial_\tau^p\partial_3^q\eta^r,_s|
+ \sum_{|p|+q=0}^{[N/2]-1}|w^{q/2}\partial_\tau^p\partial_3^q\eta_t^r,_s|<\infty\,.
\end{equation}
 Then we obtain the following a priori estimates:
 \be\label{Energy-I}
\aligned
{d \over dt} \left[  \mathcal{E}_{m,n}^{(I)}  + (1+\frac1\alpha) \mathcal{E}_{m,n}^{(II)}  \right] &\leq \Fcal_1\left(E_N^{(I)}, \,E_N^{(III)}\right) \; \text{ for }\;  |m| <N,\\
{d \over dt} \left[  \mathcal{E}_{N,0}^{(I)}  + (1+\frac1\alpha) \mathcal{E}_{N,0}^{(II)} +\mathcal{G}  \right] &\leq \Fcal_1\left(E_N^{(I)}, \,E_N^{(III)}\right)\; \text{ for }\;  |m| =N,
\endaligned
\ee
where for any $\delta>0$
\[
|\mathcal{G}|\leq \delta \mathcal{E}^{(III)}_{N,0} + C_\delta \mathcal{E}^{(III)}_{N-1,0}
\] 
as well as 
 \begin{equation}\label{e}
 \begin{split}
&E^{(III)}_N \leq \Fcal_2 \left(E_N[\eta_0, \eta_1], E_N^{(I)}, \,E_N^{(III)}, T\right), \\
\
\
&E^{(IV)}_N \leq E^{(IV)}_N[\eta_0,\eta_1] + \Fcal_3\left(E_N^{(I)}, \,E_N^{(III)}, T\right), 
\end{split}
\end{equation}
where $\Fcal_1$, $\Fcal_2$ and $\Fcal_3$ are smooth functions in their arguments. Moreover, the a priori assumption \eqref{Assumption} can be justified. 
\end{theorem}

The proof of Theorem \ref{theo1} is a direct consequence of the following two lemmas.

\begin{lemma}[Energy estimates] 
\label{energylemma}
Under the assumptions of Theorem~\ref{theo1}, one has the energy inequality \eqref{Energy-I}. 
\end{lemma} 

\begin{lemma}[Gradient and curl estimates] 
\label{energylemma1}
Under the assumptions of Theorem~\ref{theo1}, one obtains the energy bounds \eqref{e}.  
\end{lemma} 

The structure exhibited by the second-order system and the curl system \eqref{syste-0}--\eqref{recall} is fundamental in order to derive the necessary estimates. The first energy inequality  is a consequence of the wave-like structure of the second-order system, while the second energy bounds will follow from the vorticity equations. 
We observe that the energy functionals incorporate suitable powers of the weight function $w$. Higher powers are required for normal derivatives, 
while no such loss is encountered for tangential derivatives. The same algebraic structure of the change in weights has been identified for the non-relativistic flows in \cite{JM2}. 

While there is some similarity to the proof for the non-relativistic Euler flows as done in \cite{JM2}, our proof here involves some new ingredients and  the estimates are not the same. The paper \cite{JM2} derived the estimates for $\partial_t\eta$ and the full gradient of $\eta$ from the energy estimates of the second-order hyperbolic system at the expense of loosing the positivity of the curl part in the energy and the method therein compensated a lost curl energy by an auxiliary estimate from the curl equation.  The curl equation for the non-relativistic Euler flows is rather simple and elegant  (i.e.~almost an ODE) in Lagrangian coordinates, which is one of key ingredients used in \cite{JM2}, but such a simple structure does not seem to be available for the relativistic Euler equations. To get around this difficulty present for the relativistic Euler flows  even at the formal level, we obtain the estimates for $\partial_t\eta$ and the divergence of $\eta$ via the energy estimates at the expense of loosing the positivity of the full tangential derivative terms and recover the full gradient estimates from the relativistic vorticity equation. This new scheme is applied to the non-relativistic Euler equations and it gives an alternative way of deriving the estimates. 

We observe that Theorem \ref{theo1} can be extended to a larger class of quasilinear hyperbolic systems inheriting the same leading-order structure as in \eqref{syste-0} and \eqref{curl-eta-new1}, so long as the coefficient matrices and tensors satisfy suitable algebraic conditions such as symmetry, anti-symmetry, and positive definiteness.
For instance, taking into account lower-order forcing term such as a gravitational coupling or damping terms would not add any further difficulty at this level.

\subsection{Hardy inequality and embedding of weighted Sobolev spaces}

Before we derive the energy estimates, we recall the following useful Hardy inequality and embedding results. 
First of all, for the Hardy inequality we have the following \cite{KMP}. 

\begin{lemma} \label{hardy} (Hardy inequality) 
Let  $k $ be a real number and $g$ a function satisfying $\int_0^1  s^k (g^2 + g'^2) ds < \infty $.

If $k > 1$, then we have $ \int_0^1  s^{k-2}  g^2 ds  \leq C \int_0^1  s^k (g^2 + |g'|^2) ds   $. 

If $k < 1$, then $g$ has a trace at $x=0$  and 
$ \int_0^1  s^{k-2}  (g - g(0))^2 ds  \leq C \int_0^1  s^k   |g'|^2  ds$. 
 \end{lemma}

Note that using Lemma \ref{hardy} with $k = \alpha +1$, we get  
  \begin{equation}  \label{hard1}
 \int_\Omega  w^{\alpha - 1 } | v  |^2 dx \,  
  \leq     C \int_\Omega  [w^{\alpha + 1 } |\partial_3 v  |^2       +    w^{\alpha+1} |v|^2 ]    dx.  \,
\end{equation} 
We will also use the following variant of Hardy inequality: for any fixed $\delta>0$, 
  \begin{equation}  \label{hard2}
 \int_\Omega  w^{\alpha - 1 } | v  |^2 dx \,  
  \leq     \delta \int_\Omega w^{\alpha + 1 } |\partial_3 v  |^2 dx   +  C_\delta\int_\Omega     w^{\alpha+1}|v|^2    dx.  \,
\end{equation} 

The above energy functionals induce a family of weighted Sobolev spaces. It is convenient to introduce the function spaces $X^{\alpha,b}$, $Y^{\alpha,b}$, $Z^{\alpha,b}$ to discuss the embedding results: 
\begin{equation}\label{XY}
\begin{split}
X^{\alpha,b}&\equiv \{w^{\frac{\alpha}{2}}F\in L^2(\Omega) : \int_{\Omega}w^{\alpha+n}|\partial_\tau^m\partial_3^n F|^2 dx<\infty\,,\,0\leq |m|+n\leq b\},\\
Y^{\alpha,b}&\equiv \{w^{\frac{1+\alpha}{2}}D_\eta F\in L^2(\Omega) : \int_{\Omega}w^{1+\alpha+n}|D_\eta \partial_\tau^m\partial_3^n F|^2 dx<\infty\,,\,0\leq |m|+n\leq b\},\\
Z^{\alpha,b}&\equiv \{w^{\frac{1+\alpha}{2}} F\in L^2(\Omega) : \int_{\Omega}w^{1+\alpha+n} | \partial_\tau^m\partial_3^n F|^2 dx<\infty\,,\,0\leq |m|+n\leq b\}.
\end{split}
\end{equation}

Then as an application of the Hardy type embedding of weighted Sobolev spaces \cite{KMP}, we obtain the 
embedding of $X^{\alpha,b}$,  $Y^{\alpha,b}$,  $Z^{\alpha,b}$  into the standard Sobolev spaces $H^s$ for  sufficiently smooth $w$. 

\begin{lemma}\label{emb} For $b\geq \lceil\alpha\rceil$, 
\begin{equation*}
\|F\|_{H^{\frac{b-\alpha}{2}}}\precsim \|F\|_{X^{\alpha,b}}\,.
\end{equation*}
In particular, for $b\geq [\alpha]+4$, 
\begin{equation*}
\|F\|_\infty \precsim  \|F\|_{X^{\alpha,b}}\,.
\end{equation*}
We have 
the similar embeddings for $Y^{\alpha,b}$ and $Z^{\alpha,b}$: for $b\geq \lceil\alpha\rceil+1$,  
\begin{equation*}
\|DF\|_{{H}^{\frac{b-\alpha-1}{2}}}\precsim \|F\|_{Y^{\alpha,b}}\;\;\text{ and  }\;\;
\|F\|_{{H}^{\frac{b-\alpha-1}{2}}}\precsim \|F\|_{Z^{\alpha,b}}\,.
\end{equation*}
\end{lemma}

We observe that the a priori bound in \eqref{Assumption} in Theorem \ref{theo1} can be justified in our energy function spaces by using Lemma \ref{hardy} and Lemma \ref{emb}, in other words  
 $|w^{q/2}\partial_\tau^p\partial_3^q\eta^r,_s|$ and 
 $|w^{q/2}\partial_\tau^p\partial_3^q\eta_t^r,_s|$ for $0\leq |p|+q\leq [N/2]$ are bounded by $E_N$.  

The remaining part of this section is devoted to the proof of Lemma \ref{energylemma} and Lemma \ref{energylemma1}.


\subsection{Proof of Lemma \ref{energylemma}} The energy inequality \eqref{Energy-I} is due to the symmetric structure of  the reformulation \eqref{syste-0}. While there is some similarity to the proof for non-relativistic Euler  as done in \cite{JM2}, our proof here is not the same. Unlike in \cite{JM2}, we will not keep the precise curl structure at the level of the energy estimates, but aim to control the divergence part only at this point. Then the full energy will be recovered by exploiting the curl equations. We notice that the obvious difference lies in that $B^j_i$ is a symmetric positive definite matrix and $C^k_{ij}$ is a symmetric tensor for the current case, while $B^j_i=\delta^j_i$ and $C^k_{ij}=0$ for the non-relativistic Euler case.  

The proof consists of three steps: the zeroth order estimate, the derivation of high order equations, and the high order estimates. Let us start with the zeroth order estimate. 

\

\noindent\textbf{Step 1 - the zeroth order estimate:} Multiply \eqref{syste-0} by $\eta_{t}^j$ and integrate to get 
\[
\int_{\Omega} w^{\alpha}\eta_t^j B^j_i\eta_{tt}^i dx + \int_{\Omega} w^{\alpha+1} \eta_t^j C^k_{ij}\del_k\eta_t^i dx + \int_{\Omega}\eta_t^j  \del_k \big( w^{\alpha+1} A_j^k J^{-1/\alpha} \big) dx =0. 
\]
The first and second terms can be written as 
\[
\begin{split}
&\int_{\Omega} w^{\alpha}\eta_t^j B^j_i\eta_{tt}^i dx  = \frac12\frac{d}{dt}\int_{\Omega} w^{\alpha}\eta_t^j B^j_i\eta_{t}^i dx - \frac12 \int_{\Omega} w^{\alpha}\eta_t^j \del_tB^j_i\eta_{t}^i dx\\
 &\int_{\Omega} w^{\alpha+1} \eta_t^j C^k_{ij}\del_k\eta_t^i dx = -  \frac12\int_{\Omega} \eta_t^j \del_k (w^{\alpha+1} C^k_{ij})\eta_t^i dx \end{split}
\]
by using the symmetry relation $B^j_i=B^i_j$ and $C^k_{ij}=C^k_{ji}$. The third term can be written as 
\[
\int_{\Omega} \eta_t^j  \del_k \big( w^{\alpha+1} A_j^k J^{-1/\alpha} \big) dx  = \frac{d}{dt} \int_{\Omega} \alpha w^{\alpha+1}J^{-1/\alpha} dx
\]
by using \eqref{DJ}. Thus \eqref{Energy-I} is valid for $m=0$ and $n=0$ in the energy. 

\

\noindent\textbf{Step 2 - the derivation of high order equations:} Let $m$ and $n$  for $1\leq |m|+n\leq N$ be fixed. 
Taking $ \del_\tau^m \del_3^n$ of  $w^{-\alpha}\cdot\eqref{syste-0}$ and by multiplying it back by $w^{\alpha+n}$, we first obtain 
 \begin{equation}\label{3nm}
 \begin{split}
 &w^{\alpha+n}B^j_i\del_\tau^m\del_3^n\eta_{tt}^i+ \sum_{|p|+q<|m|+n} c_{p,q} w^{\alpha+n} \partial_\tau^{m-p}\partial_3^{n-q} B^j_i\partial_\tau^p\partial_3^q\eta_{tt}^i \\
 & +w^{1+\alpha+n}C^k_{ij}\del_k\del_\tau^m\del_3^n\eta_{t}^i + \sum_{|p|+q<|m|+n} c_{p,q} w^{\alpha+n}  \partial_\tau^{m-p}\partial_3^{n-q} \left[w C^k_{ij}\right]\del_k\partial_\tau^p\partial_3^q\eta_{t}^i 
 \\ 
&+\underline{\underline{w^{\alpha+n}\partial_\tau^m\partial_3^n \left(  w \del_k (A^k_j J^{-1/\alpha}) +(1+\alpha) \del_k w A^k_j J^{-1/\alpha}  \right) }}=0
\end{split}
\end{equation}

We first claim that the last double-lined term in \eqref{3nm} can be written as follows: 
\be\label{claim3}
\begin{split}
&w^{\alpha+n}\,\del_\tau^m\del_3^n\left(  w \del_k (A^k_j J^{-1/\alpha}) +(1+\alpha) \del_k w A^k_j J^{-1/\alpha}  \right) \\
&\quad=-(1+\frac1\alpha) \del_k \left( w^{1+n+\alpha}   J^{-\frac1\alpha} A^k_j \text{div}_\eta\del_\tau^m\del_3^n\eta  \right)\\
&\quad\quad+(1+\alpha)w^{\alpha+n} J^{-\frac1\alpha}\del_3w
(A^3_jA^\sigma_r-A^3_rA^\sigma_j)\del_\tau^m\del_3^n\eta^r,_\sigma +w^{\alpha+n} \mathcal{R}_{m,n}, 
\end{split}
\ee
where $\mathcal{R}_{m,n}$ consists of lower order terms: $\mathcal{R}_{m,n}=$
\be\label{Rmn}
\begin{split}
\mathcal{R}_{m,n}\Big(&\del_\tau^{m-(p+q)}\del_3^{n-(i+j)}w \del_\tau^p \del_3^iD^2\eta \del_\tau^q \del_3^jD\eta,\,  \del_\tau^{m-(p+q)}\del_3^{n-(i+j)}\del_\sigma w\del_\tau^p \del_3^{i}
D\eta \del_\tau^q\del_3^jD\eta, \\
&\del_\tau^q \del_3^j w\del_\tau^p \del_3^i D\eta \del_\tau^{m-(p+q)} \del_3^{n-(i+j)} D\del_\sigma \eta,\,\del_\tau^q \del_3^jD w\del_\tau^p \del_3^i D\eta \del_\tau^{m-(p+q)} \del_3^{n-(i+j)}D \eta; \\
&\; 0\leq|p|+ i\leq |m|+n-1 \,; 1\leq |q|+ j\leq |m|+n\,; i+j\leq n\,; p+q \leq m  \Big). 
\end{split}
\ee 

We observe that the structure encoded in \eqref{claim3} is different from the one in \cite{JM2}. A new aspect is that instead of looking at the gradient of  the full gradient plus divergence minus curl as suggested by the following identity
\[
\partial_l(A_i^k{J}^{-1/\alpha})=-{J}^{-1/\alpha}A^k_r\, [D_\eta\partial_l\eta]^i_r
-\tfrac{1}{\alpha}{J}^{-1/\alpha}A^k_i\,\text{div}_\eta\partial_l\eta-{J}^{-1/\alpha}A^k_r\,
[\text{Curl}_\eta\partial_l\eta]^r_i
\]
we will make use of the structure of  the gradient of the divergence. To make it precise, first note that 
\be\label{structure1}
 \begin{split}
 \partial_l(A_i^k{J}^{-1/\alpha}) 
&=-{J}^{-1/\alpha}A^k_rA^s_i\partial_l\eta^r,_s-\tfrac{1}{\alpha}{J}^{-1/\alpha}A^k_iA^s_r\partial_l
\eta^r,_s\\
&= -(1+\frac1\alpha){J}^{-1/\alpha}A^k_i\text{div}_\eta \del_l \eta + {J}^{-1/\alpha}\left[A^k_i A^s_r - A^k_rA^s_i\right] \del_l\eta^r,_s
 \end{split}
\ee
and moreover,  
\be\label{structure2}
 \begin{split}
 \del_l\partial_k(A_i^k{J}^{-1/\alpha})
&=-(1+\frac1\alpha) \del_k\left({J}^{-1/\alpha}A^k_i\text{div}_\eta \del_l \eta\right) \\
&\quad\quad+ \del_k\left[ {J}^{-1/\alpha}A^k_i A^s_r -  {J}^{-1/\alpha}A^k_rA^s_i\right] \del_l\eta^r,_s \\
 \end{split}
\ee
after the cancelation due to the symmetry in $k,s$: $\left[A^k_i A^s_r - A^k_rA^s_i\right] \del_k\del_l\eta^r,_s =0$. We observe that the second term in \eqref{structure2} is lower order. Based on \eqref{structure1} and \eqref{structure2}, we will establish the following equivalent expression to \eqref{claim3}: 
\be\label{claim2}
\begin{split}
&\del_\tau^m\del_3^n\left(  w \del_k (A^k_j J^{-1/\alpha}) +(1+\alpha) \del_k w A^k_j J^{-1/\alpha}  \right) \\
&\quad=-(1+\frac1\alpha) \left[ w \del_k \left(  J^{-\frac1\alpha} A^k_j \text{div}_\eta \del_\tau^m\del_3^n\eta  \right)
+(1+n+\alpha) \del_k w  J^{-\frac1\alpha} A^k_j \text{div}_\eta\del_\tau^m\del_3^n\eta
\right]
\\
&\quad\quad+(1+\alpha) J^{-\frac1\alpha}\del_3w
(A^3_jA^\sigma_r-A^3_rA^\sigma_j)\del_\tau^m\del_3^n\eta^r,_\sigma + \mathcal{R}_{m,n}
\end{split}
\ee

We will present the details for normal derivatives ($m=0$) on how the weight structure changes and move onto tangential and mixed derivatives. Our first claim is that 
\be\label{claim0}
\begin{split}
&\del_3^n\left(  w \del_k (A^k_j J^{-1/\alpha}) +(1+\alpha) \del_k w A^k_j J^{-1/\alpha}  \right) \\
&\quad=-(1+\frac1\alpha) \left[ w \del_k \left(  J^{-\frac1\alpha} A^k_j \text{div}_\eta\del_3^n\eta  \right)
+(1+n+\alpha) \del_k w  J^{-\frac1\alpha} A^k_j \text{div}_\eta\del_3^n\eta
\right]
\\
&\quad\quad+(1+\alpha) J^{-\frac1\alpha}\del_3w
(A^3_jA^\sigma_r-A^3_rA^\sigma_j)\del_3^n\eta^r,_\sigma + \mathcal{R}_{0,n}, 
\end{split}
\ee
where $\mathcal{R}_{0,n}$ consists of lower order terms: for $n\geq 1$
\be\label{Rn}
\begin{split}
\mathcal{R}_{0,n}=\mathcal{R}_{0,n}(&\del_3^{n-(i+j)}w \del_3^iD^2\eta \del_3^jD\eta,  \del_3^{n-(i+j)}\del_\sigma w \del_3^{i}D\eta \del_3^jD\eta, \\
& \del_3^j w \del_3^i D\eta \del_3^{n-(i+j)} D\del_\sigma \eta, \del_3^jD w \del_3^i D\eta  \del_3^{n-(i+j)}D \eta; \\
&\;0\leq i\leq n-1\,; 1\leq  j\leq n\,; i+j\leq n)
\end{split}
\ee

We will establish \eqref{claim0} inductively.

\underline{$\ast$ Case of $n=1$ in \eqref{claim0}.} Note that 
\be\label{l}
\begin{split}
&\del_3\left(  w \del_k (A^k_j J^{-1/\alpha}) +(1+\alpha) \del_k w A^k_j J^{-1/\alpha}  \right) \\
&\quad= w\del_3\del_k(A^k_j J^{-\frac1\alpha}) + \del_3 w \del_k (A^k_j J^{-\frac1\alpha})  + (1+\alpha) \del_k w \del_3(A^k_j J^{-\frac1\alpha} ) \\
&\quad\quad\quad+ 
(1+\alpha)\del_3 \del_k w A^k_j J^{-\frac1\alpha}
\end{split}
\ee
The second term in the right hand side of \eqref{l} is not lower order with respect to the weight. We rewrite it as 
\[
\begin{split}
 &\del_3 w \del_k (A^k_j J^{-1/\alpha})\\
 &= -(1+\frac{1}{\alpha}) J^{-1/\alpha} \del_3w A^k_jA^s_r\del_k\eta^r,_s \\
 &=  -(1+\frac{1}{\alpha}) J^{-1/\alpha} \del_kw A^k_jA^s_r\del_3\eta^r,_s +(1+\frac1\alpha)J^{-1/\alpha}  (\del_\sigma w A^\sigma_jA^s_r\del_3\eta^r,_s - \del_3 wA^\sigma_j A^s_r \del_\sigma \eta^r,_s)
\end{split}
\]
Now the second and third terms in \eqref{l} together become
\[
\begin{split}
&\del_3 w \del_k (A^k_j J^{-1/\alpha})  + (1+\alpha) \del_k w \del_3(A^k_j J^{-1/\alpha} ) \\
&=-(1+\frac{1}{\alpha}) (2+\alpha)\del_kw J^{-1/\alpha}  A^k_j \text{div}_\eta \del_3\eta + 
(1+\alpha)J^{-1/\alpha}\del_kw
(A^k_jA^s_r-A^k_rA^s_j)\del_3\eta^r,_s \\
&
\quad+(1+\frac1\alpha)J^{-1/\alpha}  (\del_\sigma w A^\sigma_jA^s_r\del_3\eta^r,_s - \del_3 wA^\sigma_j A^s_r \del_\sigma \eta^r,_s)
\end{split}
\]
but then, the second term in the right hand side when $k=3$ reduces to $\del_3w
(A^3_jA^s_r-A^3_rA^s_j)\del_3\eta^r,_s  =\del_3w
(A^3_jA^\sigma_r-A^3_rA^\sigma_j)\del_3\eta^r,_\sigma$ since when $s=3$, $A^3_jA^3_r-A^3_rA^3_j=0$. 
Hence by using \eqref{structure2} for the first term in \eqref{l}, we see that \eqref{l} can be rewritten as 
\be\label{l'}
\begin{split}
&\del_3\left(  w \del_k (A^k_j J^{-1/\alpha}) +(1+\alpha) \del_k w A^k_j J^{-1/\alpha}  \right) \\
&\quad=-(1+\frac1\alpha) \left[ w \del_k \left(   J^{-\frac1\alpha} A^k_j \text{div}_\eta\del_3\eta  \right) +  (2+\alpha)\del_kw J^{-1/\alpha}  A^k_j \text{div}_\eta \del_3\eta \right] \\
&\quad\quad+(1+\alpha)J^{-\frac1\alpha}\del_3w
(A^3_jA^\sigma_r-A^3_rA^\sigma_j)\del_3\eta^r,_\sigma +\mathcal{R}_{0,1}, 
\end{split}
\ee
where 
\[
\begin{split}
\mathcal{R}_{0,1}&:=w \del_k[ {J}^{-1/\alpha}A^k_i A^s_r -  {J}^{-1/\alpha}A^k_rA^s_i] \del_3\eta^r,_s\\
&\quad+(1+\alpha)J^{-\frac1\alpha}\del_\kappa w
(A^\kappa_jA^s_r-A^\kappa_rA^s_j)\del_3\eta^r,_s + (1+\alpha)\del_3 \del_k w A^k_j J^{-1/\alpha}  \\
& \quad+(1+\frac1\alpha)J^{-1/\alpha}  \left(\del_\sigma w A^\sigma_jA^s_r\del_3\eta^r,_s - \del_3 wA^\sigma_j A^s_r \del_\sigma \eta^r,_s\right). 
\end{split}
\]
Note that 
\[
\mathcal{R}_{0,1}= \mathcal{R}_{0,1}(w D^2\eta \del_3D\eta,  \del_\sigma w D\eta \del_3D\eta, \del_3Dw D\eta,\del_3w D\eta D \del_\sigma\eta )
\]
which consists of lower order terms. This verifies \eqref{claim0}  for $n=1$.

\underline{$\ast$ Case of  $n\geq 1$ in \eqref{claim0}.}  Suppose we have \eqref{claim0}. By taking $\del_3$ of  \eqref{claim0}, we first obtain 
\[
\begin{split}
&\del_3^{n+1}\left(  w \del_k (A^k_j J^{-1/\alpha}) +(1+\alpha) \del_k w A^k_j J^{-1/\alpha}  \right) \\
&\quad=-(1+\frac1\alpha) \left[ w \del_3\del_k \left(  J^{-\frac1\alpha} A^k_j \text{div}_\eta\del_3^{n}\eta  \right) 
+ \del_3 w \del_k \left(  J^{-\frac1\alpha} A^k_j \text{div}_\eta\del_3^{n}\eta  \right)
\right] \\
&\quad
\quad-(1+\frac1\alpha)(1+n+\alpha) \left[ \del_k w \del_3\left( J^{-\frac1\alpha} A^k_j \text{div}_\eta\del_3^{n}\eta 
\right) + \del_3\del_kwJ^{-\frac1\alpha} A^k_j \text{div}_\eta\del_3^{n}\eta   \right] 
\\
&\quad\quad+(1+\alpha)\del_3\left( J^{-\frac1\alpha}\del_3w
(A^3_jA^\sigma_r-A^3_rA^\sigma_j)\del_3^n\eta^r,_\sigma \right) +\del_3 \mathcal{R}_{0,n}
\end{split}
\]
We rewrite the first three terms in the right hand side after rearrangement as 
\[
\begin{split}
&-(1+\frac1\alpha) \left[ w \del_k \left(  J^{-\frac1\alpha} A^k_j \text{div}_\eta\del_3^{n+1}\eta  \right) 
+(2+n+\alpha)\del_k w \left( J^{-\frac1\alpha} A^k_j \text{div}_\eta\del_3^{n+1}\eta 
\right) \right] \\
&-(1+\frac1\alpha) \left[ w\del_k \left(   \del_3(J^{-\frac1\alpha} A^k_jA^s_r) \del_3^{n}\eta^r,_s  \right) 
+ \del_3 w \del_3 \left(  J^{-\frac1\alpha} A^3_j A^s_r\right) \del_3^{n}\eta^r,_s    \right] \\
&- (1+\frac1\alpha)\left[\del_3 w \del_\sigma\left(J^{-\frac1\alpha} A^\sigma_j   \text{div}_\eta\del_3^{n}\eta \right) -\del_\sigma w
J^{-\frac1\alpha} A^\sigma_j \text{div}_\eta\del_3^{n+1}\eta \right] \\
&-(1+\frac1\alpha)(1+n+\alpha)\del_k w \del_3\left( J^{-\frac1\alpha} A^k_j A^s_r\right) \del_3^{n}\eta^r,_s, 
\end{split}
\]
where the first line is the main structural expression. Thus, we see that 
\[
\begin{split}
&\del_3^{n+1}\left(  w \del_k (A^k_j J^{-1/\alpha}) +(1+\alpha) \del_k w A^k_j J^{-1/\alpha}  \right) \\
&\quad=-(1+\frac1\alpha) \left[ w \del_k \left(  J^{-\frac1\alpha} A^k_j \text{div}_\eta\del_3^{n+1}\eta  \right) 
+(2+n+\alpha)\del_k w  J^{-\frac1\alpha} A^k_j \text{div}_\eta\del_3^{n+1}\eta 
 \right] \\
&\quad\quad + (1+\alpha) J^{-\frac1\alpha}\del_3w
(A^3_jA^\sigma_r-A^3_rA^\sigma_j)\del_3^{n+1}\eta^r,_\sigma +\mathcal{R}_{0,n+1}, 
\end{split}
\]
where 
\[
\begin{split}
\mathcal{R}_{0,n+1}=
&-(1+\frac1\alpha) \left[ w\del_k \left(   \del_3\left(J^{-\frac1\alpha} A^k_jA^s_r\right) \del_3^{n}\eta^r,_s  \right) 
+ \del_3 w \del_3 \left(  J^{-\frac1\alpha} A^3_j A^s_r\right) \del_3^{n}\eta^r,_s    \right] \\
&- (1+\frac1\alpha)\left[\del_3 w \del_\sigma\left(J^{-\frac1\alpha} A^\sigma_j   \text{div}_\eta\del_3^{n}\eta \right) -\del_\sigma w
J^{-\frac1\alpha} A^\sigma_j \text{div}_\eta\del_3^{n+1}\eta \right] \\
&-(1+\frac1\alpha)(1+n+\alpha)\left[ \del_k w \del_3\left( J^{-\frac1\alpha} A^k_j A^s_r\right) \del_3^{n}\eta^r,_s + \del_3\del_kwJ^{-\frac1\alpha} A^k_j \text{div}_\eta\del_3^{n}\eta   \right]  \\
&+(1+\alpha)\del_3\left( J^{-\frac1\alpha}\del_3w
(A^3_jA^\sigma_r-A^3_rA^\sigma_j) \right) \del_3^n\eta^r,_\sigma +\del_3 \mathcal{R}_{0,n}
\end{split}
\]
which recovers \eqref{claim0} and \eqref{Rn} for $n+1$. 

We will now move onto the tangential and mixed derivatives in \eqref{claim2}. We first verify \eqref{claim2} for $n=0$. 

\underline{$\ast$ Case of $|m|\geq 1$ and $n=0$ in \eqref{claim2}.} Let us start with $|m|=1$ and $n=0$. 
\be\label{ll}
\begin{split}
&\del_\tau\left(  w \del_k (A^k_j J^{-1/\alpha}) +(1+\alpha) \del_k w A^k_j J^{-1/\alpha}  \right) \\
&\quad= w\del_\tau\del_k(A^k_j J^{-1/\alpha}) + \del_\tau w \del_k (A^k_j J^{-1/\alpha})  + (1+\alpha) \del_k w \del_\tau(A^k_j J^{-1/\alpha} ) \\
&\quad\quad\quad+ 
(1+\alpha)\del_\tau \del_k w A^k_j J^{-1/\alpha}
\end{split}
\ee
Then the second term in in the right hand side of \eqref{ll} is indeed lower order since $\del_\tau w$ behaves like $w$.   Hence we can rewrite it as 
\be\label{l=12}
\begin{split}
&\del_\tau\left(  w \del_k (A^k_j J^{-1/\alpha}) +(1+\alpha) \del_k w A^k_j J^{-1/\alpha}  \right) \\
&=-(1+\frac1\alpha) \left[ w\del_k \left(  J^{-\frac1\alpha} A^k_i \text{div}_\eta\del_\tau\eta  \right)+(1+\alpha) \del_k w  J^{-\frac1\alpha} A^k_i \text{div}_\eta\del_\tau\eta   \right] \\
&\quad+(1+\alpha)
 J^{-\frac1\alpha}\del_3w
(A^3_jA^\sigma_r-A^3_rA^\sigma_j)\del_\tau\eta^r,_\sigma +\mathcal{R}_{1,0}, 
\end{split}
\ee
where
\be
\begin{split}
\mathcal{R}_{1,0}&= w \del_k[ {J}^{-\frac1\alpha}A^k_i A^s_r -  {J}^{-\frac1\alpha}A^k_rA^s_i] \del_\tau\eta^r,_s + 
(1+\alpha) J^{-\frac1\alpha}\del_\sigma w
(A^\sigma_jA^s_r-A^\sigma_rA^s_j)\del_\tau\eta^r,_s \\
& \quad -(1+\frac{1}{\alpha}) J^{-1/\alpha} \del_\tau w A^k_jA^s_r\del_k\eta^r,_s+ (1+\alpha) \del_\tau \del_k w A^k_j J^{-1/\alpha}
\end{split}
\ee
We observe $\mathcal{R}_{1,0}$ can be put into the following form
\[
\mathcal{R}_{1,0}=\mathcal{R}_{1,0} (w D^2\eta \del_\tau D\eta, \del_\sigma w  D\eta  D^2\eta, \del_\tau Dw D\eta   )
\]
which consists of essentially lower order terms with respect to the derivatives and weights. One can take more tangential derivatives of \eqref{ll} to obtain 
\be
\begin{split}
&\del_\tau^m\left(  w \del_k (A^k_j J^{-1/\alpha}) +(1+\alpha) \del_k w A^k_j J^{-1/\alpha}  \right) \\
&=-(1+\frac1\alpha) \left[ w\del_k \left(  J^{-\frac1\alpha} A^k_i \text{div}_\eta\del_\tau^m\eta  \right)+(1+\alpha) \del_k w  J^{-\frac1\alpha} A^k_i \text{div}_\eta\del_\tau^m\eta   \right] \\
&\quad+(1+\alpha)
 J^{-\frac1\alpha}\del_3w
(A^3_jA^\sigma_r-A^3_rA^\sigma_j)\del_\tau^m\eta^r,_\sigma +\mathcal{R}_{m,0}, 
\end{split}
\ee
where $\mathcal{R}_{m,0}$ having the form in \eqref{Rmn} consists of lower order terms. 

\underline{$\ast$ Case of $|m|\geq 1$ and $n\geq 1$ in \eqref{claim2}.} The expression \eqref{claim2} can be derived by taking $\del_\tau$ consecutively of \eqref{claim0}. The point is that $\del_\tau w$ behaves like $w$, unlike the action of $\del_3$, the weight structure will not change under $\del_\tau$. Since the procedure is similar to the previous cases we omit the details. 

\

\noindent\textbf{Step 3 - High order energy estimates:} We will now perform the energy estimates for \eqref{3nm} for $1\leq |m|+n\leq N$. The energy inequality will be obtained by multiplying \eqref{3nm} by $\partial_\tau^m\partial_3^n\eta_{t}^j$ and integrating over the domain. We will derive the estimates line by line. 

{\underline{$\bullet$ The first line in \eqref{3nm}}.} The first term in  \eqref{3nm} yields the energy term corresponding $m,n$ in $E^{(I)}_N$ plus a commutator term
\[
\begin{split}
\int_{\Omega} w^{\alpha+n}\partial_\tau^m\partial_3^n\eta_t^j B^j_i\partial_\tau^m\partial_3^n\eta_{tt}^i dx  &= \frac12\frac{d}{dt}\int_{\Omega} w^{\alpha+n}\partial_\tau^m\partial_3^n\eta_t^j B^j_i\partial_\tau^m\partial_3^n\eta_{t}^i dx + R_1\\
&=\frac12\mathcal{E}^{(I)}_{m,n}+R_1, 
\end{split}
\]
where
$$
R_1=- \frac12 \int_{\Omega} w^{\alpha}\partial_\tau^m\partial_3^n\eta_t^j \del_tB^j_i\partial_\tau^m\partial_3^n\eta_{t}^i dx\precsim \int_{\Omega} w^{\alpha+n}\partial_\tau^m\partial_3^n\eta_t^j B^j_i\partial_\tau^m\partial_3^n\eta_{t}^i dx
$$
since $|\del_t B B^{-1} |$ is bounded due to the a priori bound \eqref{Assumption}. The second term in the first line of 
\eqref{3nm} yields essentially lower order nonlinear terms since $|p|+q<|m|+n$. By using \eqref{Assumption}, Lemma \ref{hardy} and Lemma \ref{emb}, one can deduce that those lower order terms are bounded by a continuous function of  $E_N^{(I)}$ and $E_N^{(III)}$.

{\underline{$\bullet$ The second line in \eqref{3nm}}.}  The first term can be written 
\[
\begin{split}
&\int_{\Omega} w^{1+\alpha+n}  \partial_\tau^m\partial_3^n\eta_{t}^j  C^k_{ij}\del_k\partial_\tau^m\partial_3^n\eta_{t}^i dx   \\
&= -\frac12 \int_{\Omega}   \partial_\tau^m\partial_3^n\eta_{t}^j \del_k \big( w^{1+\alpha+n} C^k_{ij}\big) \partial_\tau^m\partial_3^n\eta_{t}^i dx \;\text{ (by integration by parts)} \\
& \precsim \int_{\Omega} w^{\alpha+n} \del_\tau^m \del_3^n \eta_t^jB^j_i\del_\tau^m \del_3^n \eta_t ^i dx, 
\end{split}
\]
where the last step is due to \eqref{Assumption}. The second term in the second line is lower-order with respect to number of the derivatives and the weight and hence by standard nonlinear estimates using \eqref{Assumption}, Lemma \ref{hardy} and Lemma \ref{emb}, we see that it is bounded by a continuous function of $E_N^{(I)}$ and $E_N^{(III)}$.

{\underline{$\bullet$ The third line in \eqref{3nm}}.}  We will use the expression  \eqref{claim3}. Multiplying \eqref{claim3} by 
$\partial_\tau^m\partial_3^n\eta_{t}^j$ and integrating, we have 
\[
\begin{split}
&\int_{\Omega} w^{\alpha+n}\partial_\tau^m\partial_3^n\eta_{t}^j  \del_\tau^m\del_3^n\left(  w \del_k (A^k_j J^{-1/\alpha}) +(1+\alpha) \del_k w A^k_j J^{-1/\alpha}  \right)dx  \\
&\quad=-(1+\frac1\alpha)\int_{\Omega} \partial_\tau^m\partial_3^n\eta_{t}^j  \del_k \left( w^{1+n+\alpha}   J^{-\frac1\alpha} A^k_j \text{div}_\eta\del_\tau^m\del_3^n\eta  \right)dx\\
&\quad\quad+(1+\alpha)\int_{\Omega} w^{\alpha+n} \partial_\tau^m\partial_3^n\eta_{t}^j J^{-\frac1\alpha}\del_3w
(A^3_jA^\sigma_r-A^3_rA^\sigma_j)\del_\tau^m\del_3^n\eta^r,_\sigma dx \\
&\quad\quad+\int_{\Omega} w^{\alpha+n}\partial_\tau^m\partial_3^n\eta_{t}^j \mathcal{R}_{m,n} dx \\
&\quad=: (I)+(II)+(III)
\end{split}
\]
For $(I)$, by integration by parts, we obtain 
\[
\begin{split}
(I)&=(1+\frac1\alpha)\int_{\Omega} \del_k \partial_\tau^m\partial_3^n\eta_{t}^j  \left( w^{1+n+\alpha}   J^{-\frac1\alpha} A^k_j\text{div}_\eta\del_3\eta  \right)dx\\
&=\frac12\frac{d}{dt} (1+\frac1\alpha)\int_{\Omega} w^{1+n+\alpha}   J^{-\frac1\alpha} |\text{div}_\eta\partial_\tau^m\partial_3^n\eta |^2  dx \\
&\quad-\frac12 (1+\frac1\alpha)\int_{\Omega} w^{1+n+\alpha}  \del_t( J^{-\frac1\alpha}) |\text{div}_\eta\partial_\tau^m\partial_3^n\eta |^2  dx\\
&\quad- (1+\frac1\alpha)\int_{\Omega} \del_k \partial_\tau^m\partial_3^n\eta^j w^{1+n+\alpha}   J^{-\frac1\alpha} \del_tA^k_j\text{div}_\eta\partial_\tau^m\partial_3^n\eta dx
\end{split}
\]
The first term is the energy term $\mathcal{E}^{(II)}_{m,n}$ in $E_N^{(II)}$ and the last two terms are commutators.  Since $J$, $\del_tJ$, $\del_tA$ are bounded due to \eqref{Assumption}, those commutators are  bounded by 
$\int_{\Omega} w^{1+n+\alpha}|\text{div}_\eta \partial_\tau^m\partial_3^n\eta|^2dx$ and $\int_{\Omega} w^{1+n+\alpha}| D\partial_\tau^m\partial_3^n\eta |^2dx$, which are in turn bounded by $E_N^{(III)}$. For $(II)$, we divide into cases. If $n\geq 1$, since 
\[
\frac{(II)}{1+\alpha}= \int_{\Omega} w^{\frac{\alpha+n}{2}} \partial_\tau^m\partial_3^n\eta_{t}^j J^{-\frac1\alpha}\del_3w
(A^3_jA^\sigma_r-A^3_rA^\sigma_j)  w^{\frac{\alpha+n}{2}}  \del_\sigma\del_\tau^m \del_3^{n-1}\eta^r,_3 dx
\]
we deduce that it's bounded by $\int_{\Omega} w^{n+\alpha}| \partial_\tau^m\partial_3^n\eta_t|^2dx$ and $\int_{\Omega} w^{n+\alpha}|\del_\sigma\del_\tau^m \del_3^{n-1}D\eta |^2dx$, which are in turn bounded by $E_N^{(I)}$ and $E_N^{(III)}$. If $n=0$, however, it is not immediate to see that it can be controlled by our energy because it involves full tangential derivatives with only $w^\alpha$ weight. For $1\leq |m|\leq N-1$, 
\[
\begin{split}
\left| \frac{(II)}{1+\alpha}\right| &= \left|\int_{\Omega} w^{\frac{\alpha}{2}} \partial_\tau^m\eta_{t}^j J^{-\frac1\alpha}\del_3w
(A^3_jA^\sigma_r-A^3_rA^\sigma_j)  w^{\frac{\alpha}{2}}  \del_\tau^m \eta^r,_\sigma dx\right|\\
&\precsim \int_{\Omega} w^\alpha |\partial_\tau^m\eta_{t} |^2 dx + \int_{\Omega} w^\alpha | \del_\tau^{m+1} \eta|^2dx\\
&\precsim \int_{\Omega} w^\alpha |\partial_\tau^m\eta_{t} |^2 dx + \int_{\Omega} w^{\alpha +2 } | \del_\tau^{m+1}D \eta|^2dx +  \int_{\Omega} w^{\alpha +1 } | \del_\tau^{m+1}\eta|^2dx \\
&\precsim E_N^{(I)}+ E_N^{(III)} \quad (\text{since } |m|\leq N-1), 
\end{split}
\]
where we have used Hardy inequality \eqref{hard1}. Now let $n=0$ and $|m|=N$, namely full tangential derivatives.   The previous trick via Hardy inequality would not directly work for this case. We will aim to show the second inequality in \eqref{Energy-I} with the new term $\mathcal{G}$. We will write it as two terms first 
\be\label{ii}
\begin{split}
\frac{(II)}{1+\alpha}&= \int_{\Omega}   \del_\tau^m\eta_t^j w^\alpha J^{-\frac1\alpha}\del_3w
A^3_jA^\sigma_r\del_\tau^m\eta^r,_\sigma dx - \int_{\Omega}   \del_\tau^m\eta_t^j w^\alpha J^{-\frac1\alpha}\del_3w
A^3_rA^\sigma_j \del_\tau^m\eta^r,_\sigma dx\\
&=: (II)_1-(II)_2
\end{split}
\ee
and rewrite the first term $(II)_1$  by performing integration by parts in time and then space: 
\[
\begin{split}
(II)_1&=\frac{d}{dt}\int_{\Omega}   \del_\tau^m\eta^j w^\alpha \del_3wJ^{-\frac1\alpha} 
A^3_jA^\sigma_r\del_\tau^m\eta^r,_\sigma dx - \int_{\Omega}   \del_\tau^m\eta^j w^\alpha \del_3wJ^{-\frac1\alpha} 
A^3_jA^\sigma_r\del_\tau^m\eta_t^r,_\sigma dx\\
&\quad-\int_{\Omega}   \del_\tau^m\eta^j w^\alpha \del_3w\del_t(J^{-\frac1\alpha} 
A^3_jA^\sigma_r)\del_\tau^m\eta^r,_\sigma dx \\
&=\frac{d}{dt}\int_{\Omega}   \del_\tau^m\eta^j w^\alpha \del_3wJ^{-\frac1\alpha} 
A^3_jA^\sigma_r\del_\tau^m\eta^r,_\sigma dx + \boxed{\int_{\Omega}   \del_\tau^m\eta^j,_\sigma w^\alpha \del_3wJ^{-\frac1\alpha} 
A^3_jA^\sigma_r\del_\tau^m\eta_t^r dx}\\
&\quad + \int_{\Omega}   \del_\tau^m\eta^j \del_\sigma (w^\alpha \del_3wJ^{-\frac1\alpha} 
A^3_jA^\sigma_r)\del_\tau^m\eta_t^r dx-\int_{\Omega}   \del_\tau^m\eta^j w^\alpha \del_3w\del_t(J^{-\frac1\alpha} 
A^3_jA^\sigma_r)\del_\tau^m\eta^r,_\sigma dx \\
\end{split}
\]
Note that the boxed term is the same as the other term $(II)_2$  in \eqref{ii}, so they cancel out. Hence 
\[
\begin{split}
\frac{(II)}{1+\alpha}&= \frac{d}{dt}\int_{\Omega}   \del_\tau^m\eta^j w^\alpha \del_3wJ^{-\frac1\alpha} 
A^3_jA^\sigma_r\del_\tau^m\eta^r,_\sigma dx \\
& + \int_{\Omega}   \del_\tau^m\eta^j \del_\sigma (w^\alpha \del_3wJ^{-\frac1\alpha} 
A^3_jA^\sigma_r)\del_\tau^m\eta_t^r dx-\int_{\Omega}   \del_\tau^m\eta^j w^\alpha \del_3w\del_t(J^{-\frac1\alpha} 
A^3_jA^\sigma_r)\del_\tau^m\eta^r,_\sigma dx \\
&=:\frac{d}{dt} (i) + (ii) -(iii)
\end{split}
\]
We can now employ the Hardy inequality  \eqref{hard1} and \eqref{hard2} for $(i)$, $(ii)$ and $(iii)$. We will present the detail for $(i)$. 
\[
\begin{split}
|(i)|&=\left|\int_{\Omega}  w^{\frac{\alpha-1}{2}} \del_\tau^m\eta^j \del_3wJ^{-\frac1\alpha} 
A^3_jA^\sigma_r w^{\frac{\alpha+1}{2}} \del_\tau^m\eta^r,_\sigma dx\right| \\
&\leq C_\theta \int_{\Omega} w^{\alpha-1}| \del_\tau^m\eta |^2 dx + \theta \int_{\Omega} w^{\alpha+1}| \del_\tau^mD\eta |^2 dx  \text{ by Cauchy-Swartz} \\
&\leq C_\theta \delta \int_{\Omega} w^{\alpha+1}| \del_\tau^mD\eta |^2 dx + C_\delta C_\theta \int_{\Omega} w^{\alpha+1} | \del_\tau^m\eta |^2 dx
\\
& \hskip5.cm  +\theta  \int_\Omega w^{\alpha+1}| \del_\tau^mD\eta |^2 dx\text{ (by \eqref{hard2}).}
\end{split}
\]
We can choose $\theta$ and $\delta$ small if necessary. This justifies the existence and estimate of $\mathcal{G}$ in  the second inequality in \eqref{Energy-I}. Estimation of $(ii)$ and $(iii)$ follows similarly by Hardy inequality: 
\[
|(ii)|+|(iii)| \precsim E^{(I)}_N+ E^{(III)}_N. 
\]
For $(III)$ containing lower order terms, by using \eqref{Assumption}, Lemma \ref{hardy} and Lemma \ref{emb}, one can deduce that it  bounded by a continuous function of  $E_N^{(I)}$ and $E_N^{(III)}$. This concludes the proof of Lemma \ref{energylemma}. 

\subsection{Proof of Lemma \ref{energylemma1}}

Let $G=\del_\tau^m\del_3^n\eta$ be given for fixed $m$ and $n$. By taking a number of derivatives of \eqref{curl-eta-new1}, we obtain 

\begin{equation} \label{G-eq0}
 U^r_i [\text{D}_\eta  \partial_t   G]_r^j - 
   [\text{D}_\eta  \partial_t G ]^i_r U^j_r  
+  \eps^2  \Gamma^2 U_i^r   (\del_t^2 G^l \del_t\eta^r -\del_t^2 G^r \del_t\eta^l  )  U_l^j = \mathcal{T}_{m,n}, 
\end{equation} 
where 
\be\label{Tmn}
\begin{split}
&\mathcal{T}_{m,n}:= \del_\tau^m\del_3^n \left[ U_i^r  X_r^l U_l^j\right] -\sum_{|p|+q\geq 1}  \del_\tau^p\del_3^q \left[U^r_i A^s_r\right] \del_\tau^{m-p}\del_3^{n-q}\del_t \eta^j,_s \\
&-\sum_{|p|+q\geq 1}  \del_\tau^p\del_3^q \left[U^j_r A^s_r\right] \del_\tau^{m-p}\del_3^{n-q}\del_t \eta^i,_s -\sum_{|p|+q\geq 1}  \del_\tau^p\del_3^q \left[\eps^2\Gamma^2 U^r_i U^j_r\del_t\eta^r  \right] \del_\tau^{m-p}\del_3^{n-q}\del_t^2 \eta^l\\
&+\sum_{|p|+q\geq 1}  \del_\tau^p\del_3^q \left[\eps^2\Gamma^2 U^r_i U^j_r\del_t\eta^l  \right] \del_\tau^{m-p}\del_3^{n-q}\del_t^2 \eta^r
\end{split}
\ee
 In turn, we integrate in time \eqref{G-eq0} to get 
\begin{equation} \label{G-eq1}
 U^r_i [\text{D}_\eta    G]_r^j - 
   [\text{D}_\eta  G ]^i_r U^j_r  
+  \eps^2  \Gamma^2 U_i^r   (\del_t G^l \del_t\eta^r -\del_t G^r \del_t\eta^l  )  U_l^j = \mathcal{S}_{m,n}, 
\end{equation} 
where $ \mathcal{S}_{m,n}$ consists of lower order terms: 
\be\label{Smn}
\begin{split}
 \mathcal{S}_{m,n}& :=  \left( U^r_i [\text{D}_\eta    G]_r^j - 
   [\text{D}_\eta  G ]^i_r U^j_r  
+  \eps^2  \Gamma^2 U_i^r   (\del_t G^l \del_t\eta^r -\del_t G^r \del_t\eta^l  )  U_l^j \right) \Big|_{t=0}\\
&\quad+  \int_0^t  \mathcal{T}_{m,n} dt +  \int_0^t \del_t(U^r_iA^s_r )  G^j,_s  dt -  \int_0^t \del_t(U^j_rA^s_r )  G^i,_s  dt  \\
&\quad+\int_0^t  \del_t(\eps^2\Gamma^2 U^r_i U^j_r\del_t\eta^r ) \del_t G^l   dt 
- \int_0^t  \del_t(\eps^2\Gamma^2 U^r_i U^j_r\del_t\eta^l ) \del_t G^r   dt
\end{split}
\ee

We will derive the estimates for $G$ by taking the matrix scalar product of \eqref{G-eq1} with  $ w^{1+\alpha+n}[\text{D}_\eta  G]$.
The first term 
gives a control of  $ [\text{D}_\eta  G]$:  
\begin{equation}
\int_{\Omega} w^{1+\alpha+n} [\text{D}_\eta  G]_r^j  U^r_i   [\text{D}_\eta  G]_i^j dx
\end{equation} 
 The second  term 
can be integrated by parts:
\begin{equation}\begin{split} \label{part}
&\int_{\Omega} w^{1+\alpha+n} [\text{D}_\eta  G]_r^i  U_r^j  [\text{D}_\eta  G]_i^jdx = \int_{\Omega} w^{1+\alpha+n}A^s_r G^i,_s U_r^jA^k_i G^j,_k dx\\  &= - \int_{\Omega} w^{1+\alpha+n}A^s_r A^k_i G^i,_{sk} U_{r}^jG^j dx -  \int_{\Omega} (w^{1+\alpha+n}),_k A^k_iA^s_r G^i,_s  U_{r}^j G^j dx \\
&\quad-    \int_{\Omega} w^{1+\alpha+n}(A^s_rU_{r}^jA^k_i),_k G^i,_s  G^j dx 
 =:(a)+(b)+(c). 
\end{split}
\end{equation} 
For $(a)$, we integrate by parts again to get 
\[
\begin{split}
(a)= \int_{\Omega} w^{1+\alpha+n}  \text{div}_\eta  G \, U^j_r[D_\eta G]^j_r dx + 
\int_{\Omega} (w^{1+\alpha+n}),_sA^s_r A^k_i G^i,_{k}  U_{r}^jG^j dx \\+ \int_{\Omega} w^{1+\alpha+n}(A^s_r A^k_i U_{r}^j),_s G^i,_{k} G^j dx 
\end{split}
\]
and hence, 
\[
\begin{split}
&\int_{\Omega} w^{1+\alpha+n} [\text{D}_\eta  G]_r^i  U_r^j  [\text{D}_\eta  G]_i^jdx  =(a)+(b)+(c)\\
&=\int_{\Omega} w^{1+\alpha+n}  \text{div}_\eta  G \, U^j_r[D_\eta G]^j_r dx \\
&\quad+ 
\int_{\Omega} (w^{1+\alpha+n}),_sA^s_r A^k_i G^i,_{k}  U_{r}^jG^j dx  - \int_{\Omega} (w^{1+\alpha+n}),_k A^k_iA^s_r G^i,_s  U_{r}^j G^j dx \\
&\quad+\int_{\Omega} w^{1+\alpha+n}(A^s_r A^k_i U_{r}^j),_s G^i,_{k} G^j dx -    \int_{\Omega} w^{1+\alpha+n}(A^s_rU_{r}^jA^k_i),_k G^i,_s  G^j dx
\end{split}
\]
It is clear that the the last two terms are lower-order. The first term in the right-hand-side is bounded by 
\[
\aligned
& \left|  \int_{\Omega} w^{1+\alpha+n}  \text{div}_\eta  G \, U_{r}^j[D_\eta G]^j_r dx\right| 
\\
& \leq \frac18 \int_{\Omega} w^{1+\alpha+n} U_{i}^r [\text{D}_\eta  G]_r^j  [\text{D}_\eta  G]_i^j dx + C \int_{\Omega} w^{1+\alpha+n} | \text{div}_\eta  G|^2 \, dx. 
\endaligned
\]
 The middle two terms need a special attention because they may not have the right weights, for instance  the second term when $s=3$ and  the third term when $k=3$ would have stronger weight $w^{\alpha+n}$ than the desired weight $w^{1+\alpha+n}$.  This can be overcome through the Hardy inequality. Here is the estimate of the third term when $k=3$. 
\[
\begin{split}
&(1+\alpha+n)\left|\int_{\Omega} w^{\alpha+n} \del_3w A^3_i [\text{D}_\eta  G]_r^i   U_{r}^j G^j dx\right| \\
&\leq \frac{1}{16} \int_{\Omega} w^{1+\alpha+n} U_i^r [\text{D}_\eta  G]_r^j  [\text{D}_\eta  G]_i^j dx + C \int_{\Omega} w^{\alpha+n-1}|G|^2 dx
\\
&  \hskip3.cm \text{ (by the Cauchy-Schwarz inequality)}\\
&\leq\frac{1}{16} \int_{\Omega} w^{1+\alpha+n} U_i^r [\text{D}_\eta  G]_r^j  [\text{D}_\eta  G]_i^j dx +\delta \int_{\Omega} w^{1+\alpha+n}|DG|^2 dx
\\
& \quad + C_\delta \int_{\Omega} w^{1+\alpha+n}|G|^2 dx \text{ by \eqref{hard2}}\\
&\leq  \frac18 \int_{\Omega} w^{1+\alpha+n} U_i^r [\text{D}_\eta  G]_r^j  [\text{D}_\eta  G]_i^j dx +C_\delta \int_{\Omega} w^{1+\alpha+n}|G|^2 dx, 
\end{split}
\]
where the last step can be achieved by choosing $\delta>0$ appropriately.  

The last term in the left-hand-side of \eqref{G-eq1} can be bounded by 
\[
\begin{split}
&\left|  \int_{\Omega} w^{1+\alpha+n}     \eps^2  \Gamma^2 U_i^r   (\del_t G^l \del_t\eta^r -\del_t G^r \del_t\eta^l  )  U_l^j  [\text{D}_\eta  G]^j_i dx\right|\\
  & \quad\leq \frac18\int_{\Omega} w^{1+\alpha+n} U_i^r [\text{D}_\eta  G]_r^j  [\text{D}_\eta  G]_i^j dx + C \int_{\Omega} w^{1+\alpha+n} |\del_t G|^2 dx
\end{split}
\]

It now remains to estimate the right-hand-side of \eqref{G-eq1}: 
\[
\int_{\Omega} w^{1+\alpha+n}[\text{D}_\eta  G]  \mathcal{S}_{m,n}  dx \leq  \frac18\int_{\Omega} w^{1+\alpha+n} [\text{D}_\eta  G]_r^j  U^r_i   [\text{D}_\eta  G]_i^j dx + C\int_{\Omega}  w^{1+\alpha+n} |  \mathcal{S}_{m,n}|^2  dx
\]
Note that $\mathcal{S}_{m,n}$ consists of initial data and the time integral of lower order terms. Standard nonlinear estimates by using \eqref{Assumption}, Lemma \ref{hardy} and Lemma \ref{emb} and the integration by parts in time when necessary (for instance, see \cite{JM2}) yield  
\[
\int_{\Omega}  w^{1+\alpha+n} |  \mathcal{S}_{m,n}|^2  dx \leq \mathcal{F}(E_N[\eta_0,\eta_1], E_N^{(I)}, E_N^{(III)}, t), 
\]
where $\mathcal{F}$ is a smooth function. This establishes the first inequality of \eqref{e}.

We further examine the dependence of initial data on $\mathcal{S}_{m,n}$. It contains some terms depending on the initial data: $ [U^r_i [\text{D}_\eta    G]_r^j - 
   [\text{D}_\eta  G ]^i_r U^j_r  
+  \eps^2  \Gamma^2 U_i^r   (\del_t G^l \del_t\eta^r -\del_t G^r \del_t\eta^l  )  U_l^j ]\big|_{t=0}$ plus some functions of  $ \del_\tau^p\del_3^q \text{Curl}_\eta\chi |_{t=0} \cdot t  $ for $0\leq |p|\leq m$ and $0\leq q\leq n$, which come from $X$ in $\mathcal{T}_{m,n}$ -- see \eqref{Smn}, \eqref{Tmn}, and \eqref{URX}. Thus we needed the initial boundedness of not only $E_N^{(I)}$ and $E_N^{(III)}$ but also $E_N^{(IV)}$.  We observe that $E_N^{(IV)}$ contains one more time derivative than $E_N^{(III)}$ and we cannot recover it by the estimates that have been presented so far. In order to estimate $E_N^{(IV)}$, we will directly use \eqref{curl-chi}. Then since 
\[
\aligned
& \del_\tau^m\del_3^n\text{Curl}_\eta\chi 
\\
& = \del_\tau^m\del_3^n\text{Curl}_\eta\chi \big|_{t=0}+\int_0^t \del_\tau^m\del_3^n [\del_t,\text{Curl}_\eta]\chi\, ds-  \int_0^t \del_\tau^m\del_3^n\left( \Gamma^{-1}  [\text{Curl}_\eta,\Gamma]\del_t\chi\right) \,ds
\endaligned
\]
by performing integration by parts in time for the second and third terms when necessary \cite{JM2}, one can deduce that 
\[
\aligned
& \int_{\Omega} w^{1+\alpha+n} \left|\del_\tau^m\del_3^n\text{Curl}_\eta\chi\right|^2  dx \\
& \leq \int_{\Omega} w^{1+\alpha+n} \left|\del_\tau^m\del_3^n\text{Curl}_\eta\chi \big|_{t=0} \right|^2  dx + \Fcal\left(E_N^{(I)}, \,E_N^{(III)}, t\right), 
\endaligned
\]
which completes the proof the lemma. 

\section{Concluding observations} \label{sec:4}

\subsection{The Euler equations of non-relativistic fluids} \label{sec:NR}

The new a priori estimates in Theorem \ref{theo1} are trivially valid for solutions to the Euler equations of non-relativistic fluids:  
\be
\label{Euler0}
\aligned 
&\del_t \rho + \del_k (\rho \, u_k ) = 0,
\\
& \del_t ( \rho \, u_j ) + \del_k \big( \rho \, u_j u_k + p \, \delta_{jk} \big) = 0.  
\endaligned  
\ee
Note that the second--order formulation above is simplified drastically when $\eps = 0$: we find $C_{ij}^k|_{\eps = 0} \equiv 0$ and 
$$
w^\alpha \, B_i^j|_{\eps = 0} \del_t^2 \eta^i + \del_k \big( w^{1+\alpha} A_j^k J^{-1/a} \big) =0, 
$$
with 
$$ 
B_i^j |_{\eps = 0} := \delta_i^j  \Gamma^{\gamma+1} = \delta_i^j,  
$$
which leads us to the  {\sl second--order formulation in Lagrangian coordinates} for non-relativistic fluids 
\be
\label{3012}
w^\alpha \, \del_t^2 \eta^j + \del_k \big( w^{1+\alpha} A_j^k J^{-1/\alpha} \big) =0. 
\ee
Similarly, the curl equation \eqref{curl-chi} reduces to the non-relativistic curl equation when  $\eps = 0$
\be\label{curl_NR}
\text{Curl}_\eta\del_t\eta =\text{Curl} \, u_0+\int_0^t  [\del_t,\text{Curl}_\eta]\del_t\eta ds. 
\ee

The non-relativistic fluids enjoy much elegant structure as it can be seen from \eqref{3012} and \eqref{curl_NR}. We observe that based on the new estimates obtained in Lemma \ref{energylemma} and \ref{energylemma1} (of course the proof for the non-relativistic case is much simpler), one can establish the existence of the solutions to \eqref{3012} justifying Theorem \ref{theo1} corresponding to $\eps=0$ by a duality argument similar to \cite{JM2}. 

\subsection{The non-relativistic limit $\eps\to 0$}

Theorem \ref{theo1} is valid for any fixed number $\epsilon\geq 0$ and it covers both relativistic and non-relativistic fluids. The non-relativistic Euler equations are recovered by letting formally $\eps\rightarrow0$ in the relativistic Euler equations and hence, a natural question arises: can one establish the convergence of the solutions of the relativistic Euler equations indexed by $\eps$ to the solutions of the non-relativistic Euler equations when $\eps\to 0$ in the presence of vacuum? The estimates in Theorem \ref{theo1} have a uniform-in-$\eps$ bound for all  sufficiently small $\eps$, and they allure the validity of the non-relativistic limit $\eps\rightarrow 0$ at least at the formal level. A rigorous justification, of course, requires an existence theory for the relativistic Euler equations.

\subsection{Final remark} 
\label{sec:5}

As presented in the previous sections, the relativistic Euler equations exhibit  an intriguing structure and it is highly non-trivial to establish the existence of the solutions satisfying the a priori estimates given in Theorem \ref{theo1}.
 In the case where the curl becomes trivial,  there is no need to keep track of the evolution of the curl and the control of the divergence energy would suffice both for getting the estimates and for the existence theory. In that situation, the existence result  follows from our a priori estimates by a similar argument as done in \cite{JM1,JM2}.  Those cases cover, for instance, 1+1 dimensional flows and 1+3 spherically symmetric flows. However, the existence question for the general relativistic fluids in vacuum still remains open and we will leave it for future study.


\

\textbf{Acknowledgements.} The first author (JJ) was supported in part by NSF grants DMS-1212142 and DMS-1351898. 
This work was done when the second author (PLF) enjoyed the hospitality of the Courant Institute of Mathematical Sciences, New York University. PLF also acknowledge financial support from the ANR grant SIMI-1-003-01 andthe European grant ITN-642768. The third author (NM) was partially supported by the NSF grant DMS-1211806.

\

\end{document}